\documentclass{ifacconf}
     
\usepackage{graphicx,color} % include this line if your document contains figures
\usepackage{amsmath}
\usepackage{flexisym}
\usepackage{breqn}
\usepackage{amssymb}
\usepackage{mathrsfs}
\usepackage{mathtools}
\usepackage{natbib}        % required for bibliography
\usepackage{subcaption}
\usepackage{comment}
\usepackage{array}
\usepackage[justification=centering]{caption}

\newcolumntype{P}[1]{>{\centering\arraybackslash}p{#1}}
\newcolumntype{M}[1]{>{\centering\arraybackslash}m{#1}}

\usepackage{mathtools}
\mathtoolsset{showonlyrefs=true}

\newtheorem{assumption}{Assumption}
\newtheorem{problem}{Problem}

\newcommand{\reals}{\mathbb{R}}
\newcommand{\sphere}{\mathbb{S}}

\newcommand{\R}{\mathbb{R}}
\newcommand{\ones}{\textbf{1}}

\newcommand{\U}{\mathcal{U}}
\newcommand{\N}{\mathbb{N}}

\newcommand{\angularvel}[3]{\omega_{\mathcal{#1}\mathcal{#2}}^{\mathcal{#3}}}
\newcommand{\dotangularvel}[3]{\dot{\omega}_{\mathcal{#1}\mathcal{#2}}^{\mathcal{#3}}}
\newcommand{\angularvelerr}[3]{\delta\omega_{\mathcal{#1}\mathcal{#2}}^{\mathcal{#3}}}

\newcommand{\angularveldes}[3]{\tilde{\omega}_{\mathcal{#1}\mathcal{#2}}^{\mathcal{#3}}}

\newcommand{\angularaccel}[3]{\alpha_{\mathcal{#1}\mathcal{#2}}^{\mathcal{#3}}}

\newcommand{\inertia}[1]{J^{\mathcal{#1}}}

\newcommand{\inertiainv}[1]{K^{\mathcal{#1}}}

\newcommand{\specorg}[1]{\relax\ifmmode 
{\mathbf{SO}}(#1)
\else $\mathbf{SO}(#1)$
\fi}
\newcommand{\speceul}[1]{\relax\ifmmode 
{\mathbf{SE}}(#1)
\else $\mathbf{SE}(#1)$
\fi}

\newcommand{\so}[1]{\relax\ifmmode 
{\mathfrak{so}}(#1)
\else $\mathfrak{so}(#1)$
\fi}

\newcommand{\se}[1]{\relax\ifmmode 
{\mathfrak{se}}(#1)
\else $\mathfrak{se}(#1)$
\fi}

\newcommand{\skewmat}[1]{\relax\ifmmode 
{\mathcal{S}}(#1)
\else $\mathcal{S}(#1)$
\fi}

\newcommand{\skewop}[1]{\relax\ifmmode 
[#1]^{\times}
\else $[#1]^{\times}$
\fi}

\newcommand{\coordframe}[1]{\mathcal{#1}}
\newcommand{\identity}[1]{I_{#1}}
\newcommand{\rotmat}[2]{\relax\ifmmode 
{R_{\mathcal{#1}}^{\mathcal{#2}}}
\else $R_{\mathcal{#1}}^{\mathcal{#2}}$
\fi}

\newcommand{\multiGauss}[2]{\mathcal{N}(#1, #2)}

\newcommand{\dotrotmat}[2]{\dot{R}_{\mathcal{#1}}^{\mathcal{#2}}}
\newcommand{\quat}[2]{q_{\mathcal{#1}}^{\mathcal{#2}}}

\newcommand{\quatscalar}[2]{\eta_{\mathcal{#1}}^{\mathcal{#2}}}
\newcommand{\quatvector}[2]{\rho_{\mathcal{#1}}^{\mathcal{#2}}}
\newcommand{\dotquatscalar}[2]{\dot{\eta}_{\mathcal{#1}}^{\mathcal{#2}}}
\newcommand{\dotquatvector}[2]{\dot{\rho}_{\mathcal{#1}}^{\mathcal{#2}}}
\newcommand{\quaterr}[2]{\delta q_{\mathcal{#1}}^{\mathcal{#2}}}

\newcommand{\quatscalarerr}[2]{\delta\eta_{\mathcal{#1}}^{\mathcal{#2}}}
\newcommand{\quatvectorerr}[2]{\delta\rho_{\mathcal{#1}}^{\mathcal{#2}}}
\newcommand{\dotquatscalarerr}[2]{\delta\dot{\eta}_{\mathcal{#1}}^{\mathcal{#2}}}
\newcommand{\dotquatvectorerr}[2]{\delta\dot{\rho}_{\mathcal{#1}}^{\mathcal{#2}}}
\newcommand{\quatdes}[2]{\tilde{q}_{\mathcal{#1}}^{\mathcal{#2}}}

\newcommand{\norm}[1]{||#1||}

\begin{document}
\begin{frontmatter}

\title{Autonomous Satellite Rendezvous and Proximity Operations 
with Time-Constrained Sub-Optimal Model Predictive Control} 

\thanks[footnoteinfo]{Approved for public release; distribution is unlimited. Public Affairs approval \#AFRL-2022-5517.}
\thanks{GB and MH were supported by AFOSR under grant FA9550-19-1-0169, 
ONR under grants N00014-19-1-2543, N00014-21-1-2495, and N00014-22-1-2435,
and AFRL under grants FA8651-22-F-1052 and FA8651-08-D-0108 TO48.}

\author[First]{Gabriel Behrendt} 
\author[Second]{Alexander Soderlund} 
\author[Third]{Matthew Hale}
\author[Fourth]{Sean Phillips}

\address[First]{University of Florida, 
   Gainesville, FL 32611 USA (e-mail: gbehrendt@ufl.edu).}
\address[Second]{Space Vehicles Directorate, Air Force Research Laboratory, Kirtland AFB, NM 87108 USA (e-mail: alexander.soderlund.1@spaceforce.mil).}
\address[Third]{University of Florida, 
   Gainesville, FL 32611 USA (e-mail: matthewhale@ufl.edu).}
\address[Fourth]{Space Vehicles Directorate, Air Force Research Laboratory, Kirtland AFB, NM 87108 USA (e-mail: sean.phillips.9@spaceforce.mil).}

\begin{abstract}                % Abstract of not more than 250 words.
This paper presents a time-constrained model predictive control strategy for the 6 degree-of-freedom (6DOF) autonomous rendezvous and docking problem between a controllable ``deputy" spacecraft and an uncontrollable ``chief" spacecraft. The control strategy accounts for computational time constraints due to limited onboard processing speed. The translational dynamics model is derived from the Clohessy-Wiltshire equations and the angular dynamics are modeled on gas jet actuation about the deputy's center of mass. Simulation results are shown to achieve the docking configuration under computational time constraints by limiting the number of allowed algorithm iterations when
computing each input. Specifically, 
we show that upwards of~$90$\% of computations can be eliminated
from a model predictive control implementation 
without significantly harming control performance. 
%Then, to show robustness of the control strategy, we simulate the 6DOF ARPOD problem under perturbations and demonstrate that we can achieve the docking state under perturbations, even with significantly constrained computation times. 
\end{abstract}


\end{frontmatter}
%===============================================================================

\section{Introduction}

Autonomy has become increasingly necessary for space applications due to the growing number of deployed satellites and limited ability to use 
human-in-the-loop control \citep{ESAreport}. The autonomous rendezvous and docking (ARPOD) problem is one of these applications which has gained significant interest because of recent space missions such as the NASA DART mission~\citep{cheng2018} and its uses in applications such as satellite inspection~\citep{bernhard2020} and servicing~\citep{ogilvie2008}. In general, the ARPOD problem requires a deputy spacecraft to maneuver to a desired state relative to a chief spacecraft.

There are a few challenges that must be addressed in the ARPOD problem to design a 
reliable and practically-implementable control strategy. One challenge
is that the dynamics 
of the satellites are nonlinear when both the translational and attitude dynamics are considered together. 
%This nonlinearity arises from the inherent coupling between the orientation and thrust direction of the deputy spacecraft and precludes the use of classical linear control methods. 
Another challenge is that these models are subject to disturbances such as sensor noise, model mismatch, and others unique to operations in space, such as J2 perturbations and solar pressure. Furthermore, in the design of a control strategy, we must consider that input commands must be computed using space-grade hardware. Since space-grade processors must be radiation-hardened (which requires lengthy construction cycles) to handle the electronics-hostile space environment \citep{bourdarie2008near} and the design of space missions are often multi-year endeavors, there is a sizable processing capabilities gap between state-of-the-art, commercial off-the-shelf CPUs and their space-grade counterparts \citep{lovelly2017comparative, lovelly2017comparative_dissertation}. A challenge of this computational gap is that space-based missions are becoming increasingly complex, 
%as autonomous control protocols require time-sensitive observations and reactions within the space environment, 
particularly in the congested low-Earth orbit domain. Thus, we require a control strategy that (i) can handle nonlinearities in dynamics, (ii) is robust to perturbations, and (iii) can be computed by a computationally-limited satellite. 

Often, autonomous control problems are solved via state feedback methods. Model predictive control (MPC) has been widely utilized \citep{falcone2007} and is a natural fit for the ARPOD problem due to its 
robustness and 
ability to handle nonlinearities. 
In conventional MPC, an optimization problem is solved at each time step to completion, i.e., 
a new control input is computed by solving an optimization problem
until a stopping condition is reached. However, given the limited available computational speed in space, we cannot reliably use conventional MPC 
new inputs to the system may be needed before the
computations to find those inputs can be completed. 

Therefore, in this work, we use a method that we refer to as \emph{time-constrained MPC}. We consider a setting in which the underlying optimization algorithm is only allowed enough time to complete a limited number of iterations when computing each input. This constraint only allows the algorithm to make some progress toward the optimum, and, when the time constraint is reached, a sub-optimal input is applied to the system by using the optimization algorithm's most recent iterate. 

It is known in the MPC literature that the optimization problem does not need to be solved to completion in order to ensure stability of the solution~\citep{mayne2014,graichen2012, graichen2010, pavlov2019}.
%Moreover, leveraging certain problem structure and considering a particular optimization algorithm it is possible to derive a lower bound on the minimum number of optimization steps required to ensure system stability
%\citep{graichen2012, graichen2010, pavlov2019}. 
These works suggest that time-constrained MPC is a viable solution to account for nonlinear dynamics, perturbations, and limited onboard computational speed seen in the ARPOD problem, and this paper formalizes and confirms this point. 

\subsection{Summary of Contributions}
To summarize, the contributions of this paper are:
\begin{enumerate}
    \item Modeling of the ARPOD problem considering both translational and attitude dynamics  
    \item Empirical validation that time-constrained MPC 
    succeeds in docking for the 6 DoF ARPOD problem
    \item Empirical validation that time-constrained MPC 
    succeeds in docking  
    for the fully actuated 6 DoF ARPOD problem with perturbations
    \item Comparison of performance under different computation constraints
    and validation that more than~$90$\% of computations can
    be eliminated in the ARPOD problem without significantly harming performance
\end{enumerate}
\subsection{Related Work}
The ARPOD problem has been studied in a variety of settings. Some works have considered only controlling translational \citep{jewison2017guidance} or attitudinal dynamics~\citep{leomanni2014,trivailo2009}. Typically, ARPOD studies consider translational control inputs that are uncoupled from the attitudinal dynamics, such as in \citep{hogan2014attitude}. While recent studies \citep{soderlund2021autonomous, soderlund2022switching} have considered coupled translational-attitudinal systems, only two-dimensional motion between the deputy and chief was considered. Other works have developed control strategies for the more complex three-dimensional 6DOF ARPOD problem with coupled translational and attitude dynamics. These works include approaches from nonlinear control such as back-stepping \citep{sun2015, wang2019}, sliding mode~\citep{yang2019, zhou2020},  learning-based control~\citep{hu2021},
and artificial potential functions~\citep{dong2018}. 
This paper differs from those earlier works in that we do not use a static control law with user-defined gains. Instead, we use MPC to solve for our next control input. 

Model Predictive Control has been used for the 3 degree-of-freedom ARPOD problem while considering various constraints, such as collision avoidance and line-of-sight constraints~\citep{petersen2014,li2017,wang2018,richards2003,di2012}.
Considering only the translational dynamics can be convenient because, under mild assumptions, the dynamics are linear. \citep{lee2017} also considers a dual-quaternion formulation of the 6DOF ARPOD problem, and they use a piecewise affine approximation of the nonlinear system model so that the MPC problem can be solved quickly 
onboard~\citep{hartley2015}. 
Other computational considerations include using a pre-computed lookup table~\citep{malyuta2021}, or designing custom predictive controller hardware~\citep{hartley2013}. 
In this paper, we differ because we
consider a nonlinear dynamic model without any approximation and we 
explicitly address computational time constraints. 

The rest of the paper is organized as follows. Section~\ref{sec:prelim} presents preliminaries.
Then Section~\ref{sec:problem} gives a time-constrained MPC problem statement. 
%In Section~\ref{sec:dynamics} we derive the nonlinear dynamic equations for the~$6$ degree of freedom ARPOD problem. 
After that, Section~\ref{sec:results} provides simulation results, and Section~\ref{sec:conclusion} concludes.

%============================================================================================
\section{Preliminaries} \label{sec:prelim}
The Euclidean $n$-dimensional space is denoted by $\reals^{n}$. The vector $v \in \reals^{n}$ is defined as $v := (v_1, \dots, v_n)^{\top}$, where the superscript $\top$ denotes the transpose operation. Unless otherwise specified, all vectors used in this paper are \textit{physical} vectors, meaning that they exist irrespective of any coordinate frame used. For instance, the vector $v^{\mathcal{A}}$ denotes that the vector is given by the coordinates defined by the frame $\mathcal{A}$, while $v^{\mathcal{B}}$ denotes that the vector is given by the coordinates defined by the frame $\mathcal{B}$. For~$v\in\R^3$,~$\skewop{v}$ denotes the skew-symmetric operator, i.e.,
\begin{equation}
    \skewop{v} = \begin{bmatrix}
                0 & -v_3 & v_2 \\
                v_3 & 0 & -v_1 \\
                -v_2 & v_1 & 0
                \end{bmatrix}.
\end{equation}
 The $n \times n$ identity matrix is denoted by $\identity{n}$, while $\textbf{0}_{n}$ denotes the zero column vector with dimension~$n\times 1$. 
  Similarly, the~$n \times 1$-dimensional column 
  vector of ones is denoted by~$\ones_n$. The $n$-dimensional multivariate normal distribution is denoted $\multiGauss{\varphi}{\sigma}$, where $\varphi \in \reals^n$ is the mean and $\sigma \in \reals^{n \times n}$ is the covariance matrix.

%\subsection{Rotation Parameterizations and Angular Kinematics}
%\label{sec:rotations}
The Lie Group \specorg{3} is the set of all real invertible $3 \times 3$ matrices that are orthogonal with determinant $1$, i.e.,
$$\specorg{3} := \{\rotmat{}{} \in \reals^{3 \times 3} \ | \ \rotmat{}{}^{\top}\rotmat{}{} = \identity{3}, \det[\rotmat{}{}] = 1\}. $$
In this work, elements of $\specorg{3}$ will be referred to as \textit{rotation matrices}. The coordinate transformation from frame $\coordframe{B}$ to frame $\coordframe{A}$ is given by
%\begin{equation} \label{eq:coord_trans}
$v^{\mathcal{A}} = \rotmat{\mathcal{B}}{\mathcal{A}} v^{\mathcal{B}}$.
%\end{equation}
Any rotation matrix $\rotmat{}{}$ can be parameterized through the unit quaternion vector
$\quat{}{} := \begin{pmatrix} \quatscalar{}{} \\ \quatvector{}{} \end{pmatrix},$
where $\quatvector{}{} = (\rho_1, \ \rho_2, \ \rho_3)^{\top} \in \reals^3$ is the vector part and $\quatscalar{}{} \in \R$ is the scalar part. The unit quaternion is contained in the unit hypersphere, $\sphere^3$, defined with four coordinates $s := (s_1 \ s_2 \ s_3 \ s_4)^{\top} \in \reals^4$ as $\sphere^3 := \{s \in \reals^4 : s^{\top}s = 1\}.$
The inverse quaternion is 
$q^{-1} = \begin{pmatrix} \quatscalar{}{} \ {-\rho}^{\top}\end{pmatrix}^{\top}$
and the identity rotation is 
$q^{I} = \begin{pmatrix} 1 \ \mathbf{0}_3^{\top} \end{pmatrix}^{\top}.$

\begin{comment}
Multiplication of two unit quaternions $q_a = (\eta_a \ \rho^{\top}_a)^{\top}$ and $q_b = (\eta_b \ \rho^{\top}_b)^{\top}$ is carried out as
\begin{equation}
\label{eq:quaternion_multiply}
q_a \otimes q_b = \begin{pmatrix}
\eta_a \eta_b - \rho^{\top}_a \rho_b \\
\eta_a \rho_b + \eta_b \rho_a + \skewop{\rho_a} \rho_b
\end{pmatrix}.
\end{equation}
One can show that $q \otimes q^{-1} = q^{-1} \otimes q = q^{I}$. 
\end{comment}

We use $\quat{B}{A}$ to denote the quaternion corresponding to the rotation matrix $\rotmat{B}{A}$. The mapping from a quaternion $\quat{}{} = (\eta \ \rho^\top)^{\top}$ to its rotation matrix $\rotmat{}{}$ is
\begin{equation}
\label{eq:quat2rot}
\identity{3} - 2\quatscalar{}{}\skewop{\quatvector{}{}} + 2\skewop{\quatvector{}{}}\skewop{\quatvector{}{}}.
\end{equation}
We note that a quaternion $\quat{}{} = (\eta \ \rho^\top)^{\top}$ and its negative $-\quat{}{}$ represent the same physical rotation\footnote{This work has opted to adopt the convention of \citep{kuipers1999quaternions} where the passive rotation operator on a vector $v$ is represented as $q^{-1} \otimes v \otimes q$ where $v$ is a ``pure" quaternion with $0$ scalar part. This choice affects the signs of Eq.\eqref{eq:quat2rot} and the angular velocity kinematics.}.
%\subsection{Angular Velocity and Rotation Kinematics}
%\label{sec:angularvelocity}

This work denotes $\angularvel{\coordframe{E}}{\coordframe{A}}{\coordframe{B}} \in \reals^{3}$ as the angular velocity of frame $\coordframe{A}$ \textit{relative} to frame $\coordframe{E}$ but with vector components represented in frame $\coordframe{B}$. The time rate change of any rotation matrix $\rotmat{\mathcal{B}}{\mathcal{A}}$ is given as
\begin{align}
\label{eq:rotation_velocity_a}
\dotrotmat{\mathcal{B}}{\mathcal{A}} &= \rotmat{\mathcal{B}}{\mathcal{A}}\skewop{\angularvel{\coordframe{A}}{\coordframe{B}}{\coordframe{B}}} \\
\label{eq:rotation_velocity_b}
\dotrotmat{\mathcal{B}}{\mathcal{A}} &= \skewop{\angularvel{\coordframe{A}}{\coordframe{B}}{\coordframe{A}}}\rotmat{\mathcal{B}}{\mathcal{A}}.
\end{align}
\begin{comment}
The time derivative of $\angularvel{\coordframe{E}}{\coordframe{A}}{\coordframe{B}} \in \reals^{3}$ is denoted $\angularaccel{\coordframe{E}}{\coordframe{A}}{\coordframe{B}} \in \reals^{3}$. 
The time rate change of any rotation matrix $\rotmat{\mathcal{B}}{\mathcal{A}}$ is given as
\begin{align}
\label{eq:rotation_velocity_a}
\dotrotmat{\mathcal{B}}{\mathcal{A}} &= \rotmat{\mathcal{B}}{\mathcal{A}}\skewop{\angularvel{\coordframe{A}}{\coordframe{B}}{\coordframe{B}}} \\
\label{eq:rotation_velocity_b}
\dotrotmat{\mathcal{B}}{\mathcal{A}} &= \skewop{\angularvel{\coordframe{A}}{\coordframe{B}}{\coordframe{A}}}\rotmat{\mathcal{B}}{\mathcal{A}}.
\end{align}
Through a composite series of rotations, there is a relative angular velocity relation. For instance, the angular velocities over the course of four frames are related as
\begin{equation}
\label{eq:ang_vel_summation}    
\angularvel{\coordframe{A}}{\coordframe{D}}{\coordframe{A}} = \angularvel{\coordframe{A}}{\coordframe{B}}{\coordframe{A}} + \angularvel{\coordframe{B}}{\coordframe{C}}{\coordframe{A}} + \angularvel{\coordframe{C}}{\coordframe{D}}{\coordframe{A}}.
\end{equation}
Additionally, we often use the property
$$\angularvel{\coordframe{A}}{\coordframe{B}}{\coordframe{A}} = -\angularvel{\coordframe{B}}{\coordframe{A}}{\coordframe{A}}.$$
\end{comment}
Analogously, the kinematics of a quaternion $\quat{A}{B} = (\quatscalar{}{} \ \rho^{\top})^{\top}$ are related to the angular velocity as
\begin{align}
\label{eq:quaternion_velocity_c}
\begin{pmatrix} \dotquatscalar{}{} \\ \dotquatvector{}{} \end{pmatrix} &= 
-\frac{1}{2}\begin{pmatrix}
-\rho^{\top} \\ \eta\identity{3} + \skewop{\rho}
\end{pmatrix} \angularvel{\coordframe{A}}{\coordframe{B}}{\coordframe{A}}.
\end{align}

%====================================================================================================================
\section{Problem Statement} \label{sec:problem}

\subsection{Conventional Model Predictive Control} \label{sec:MPC}
In this section, we introduce conventional MPC, then explain how time-constrained MPC differs from it. 
Consider the discrete-time dynamics
%\begin{equation} \label{eq:discreteDyn}
    $x(k+1) = f_d(x(k),u(k))$,
%\end{equation}
where~$x(k) \in \chi \subset \reals^{n}$,~$u(k) \in \U \subset \R^m$, and~$f_d: \chi \times \U \rightarrow \chi$.
The goal of conventional MPC is to solve an optimal control problem at each timestep~$k$ over a finite prediction horizon of length~$N\in\mathbb{N}$, using the current state as the initial state. This process generates
the control sequence~$\textbf{u}^*(k) = \{u^*(k),u^*(k+1),\dots,u^*(N-1)\}$
and the state 
sequence~$\textbf{x}^*(k) = \{x^*(k),x^*(k+1),\dots,x^*(N)\}$. 

Then, the first input in this sequence,~$u^*(k)$, is applied to the system and this process is repeated until the end of the time horizon is reached. 
Formally, the conventional MPC problem that is solved at each time~$k$ is given
next. 
\begin{problem}[Conventional MPC] \label{prob:MPC}
\begin{small}
\begin{equation} 
    \begin{split}
        \underset{\textbf{u}(k)}{\textnormal{minimize}} & \ J(x(k),\textbf{u}(k)) = \sum^{N-1}_{i=k} \ell (x(i),u(i)) \\
        \textnormal{subject to} &\ x(k) = x(0) \\
                                &\ x(N) \in \chi_f \\
                                &\  x(i+1) = f_d(x(i),u(i)), \quad i=k,k+1,\dots,N-1 \\
                                &\ u_\min \leq u(i) \leq u_\max, \quad i=k,k+1,\dots,N-1\\
                                &\ \textbf{u}(k)\in \U_N,
    \end{split}
\end{equation}
\end{small}
where~$\ell : \chi \times \U \rightarrow \reals$
is the cost functional,~$\chi_f \subseteq \chi$ is the terminal constraint set, $u_\min$, $u_\max \in \U$ are input limits,~$\textbf{u}(k)=\{u(k),\dots,u(N-1) \}$, and~$\U_N$ is the set of admissible inputs\footnote{We refer the reader to~\citep[Definition 3.9]{grune2017} for a formal definition of an admissible control sequence.}. 
%\mh{Make sure that~$\textbf{u}(k)$,~$f_d$,~$u_{min}$, and~$u_{max}$
%are explicitly defined here.}
%\gb{Done. $f_d$ is defined under Equation~\eqref{eq:discreteDyn}}
\end{problem}

%The process of solving Problem~\ref{prob:MPC} is shown in Figure~\ref{fig:Loop}.
%
%\begin{figure}[ht]
%\centering
%\includegraphics[width=0.4\textwidth]{Figures/MPCloop.png}
%\caption{A conventional MPC loop. The continuous-time dynamics~$\dot{x}$ determine the initial state~$x(0)$ that is fed into the Receding-Horizon Controller. 
%The resulting
%optimization problem is solved, giving~$\textbf{u}^*(k)$, the first element of which,~$u^*(k)$, is applied to the system. This process repeats until a termination condition is met.} 
%\label{fig:Loop}
%\end{figure}

\subsection{Time-Constrained MPC} \label{sec:tcMPC}
In conventional MPC, Problem~\ref{prob:MPC} is iteratively solved to completion at each time~$k$ by reaching a stopping condition based on the optimality of the solution.
We use~$j_k$ to denote the number
of iterations that an optimization algorithm must complete to reach
a stopping condition in the computation of~$\textbf{u}(k)$.
In computationally constrained settings, it cannot be guaranteed that there
is time to execute all~$j_k$ desired iterations because a system input
may be needed before those computations are completed. 
%Moreover, onboard hardware must compute~$j_k$ iterates before reaching completion that can be different for each~$k$, e.g., a gradient descent algorithm must compute~$j_k$ gradient iterates before reaching the stopping condition. However, provided limited onboard computational speed we cannot guarantee that there is sufficient time to solve this problem to completion. 
Therefore, we employ the following time-constrained MPC problem 
%in which we solve the following problem
that has an explicit constraint on~$j_k$. 
\begin{problem}[Time-Constrained MPC] \label{prob:tcMPC}
\begin{small}
\begin{equation} 
    \begin{split}
        \underset{\textbf{u}(k)}{\textnormal{minimize}} & \ J(x(k),\textbf{u}(k)) = \sum^{N-1}_{i=k} \ell (x(i),u(i)) \\
        \textnormal{subject to} &\ x(k) = x(0) \\
                                &\ x(N) \in \chi_f \\
                                &\  x(i+1) = f_d(x(i),u(i)), \quad i=k,k+1,\dots,N-1\\
                                &\ u_\min \leq u(i) \leq u_\max, \quad i=k,k+1,\dots,N-1\\
                                &\ \textbf{u}(k)\in \U_N \\
                                & \ j_k \leq j_\max,
    \end{split}
\end{equation}
\end{small}
where~$j_k$ is the number of iterations computed at time~$k$ and~$j_\max$ is the maximum allowable iterations for all~$k$. 
%\mh{Same as above: Make sure every symbol is defined. Should the equation
%with~$x(i+1) = \cdots$ end with~$N-1$ instead of~$N$?}
%\gb{Can I use the same symbols I defined after Problem 1 without redefining them?}
\end{problem}

The constraint~$j_k \leq j_{max}$ models scenarios with limited onboard computational speed in which there is only enough time to complete at most~$j_\max$ iterations. 
A solution to Problem~\ref{prob:tcMPC} at time~$k$ will typically
result in a sub-optimal input and state sequences, i.e.,
%which we denote by~$\tilde{x}$ and~$\tilde{u}$, i.e., 
\begin{equation}
    \begin{split}
        \tilde{\textbf{u}}(k) &= \{\tilde{u}(k),\tilde{u}(k+1),\dots,\tilde{u}(N-1)\} \\
        \tilde{\textbf{x}}(k) &= \{\tilde{x}(k),\tilde{x}(k+1),\dots,\tilde{x}(N)\}, \\
    \end{split}
\end{equation}
respectively. Then, we apply the first input of the resulting input 
sequence, namely~$\tilde{u}(k)$, and repeat this process until the 
end of the time horizon is reached. 
This setup is different from conventional MPC in that we apply a potentially sub-optimal input due to the iteration constraint~$j_k \leq j_{max}$ in Problem~\ref{prob:tcMPC}.
Next, we derive the dynamics for the 6DOF ARPOD problem
that will be used in our setup and solution to Problem~\ref{prob:tcMPC}. 

\begin{comment}
    \begin{defn} \label{def:admissibility}
    (Admissibility): For~$N\in \N$ and~$x(0) \in \chi$, a control sequence, denoted by~$\textbf{u}(k)=\{u(k),u(k+1),\dots,u(N-1)\} \in \U_N$, is said to be admissible for state~$x(0)$ if
    \begin{enumerate}
        \item $u(k) \in \U$ for all $k \in \{0,1,\dots,N\}$
        \item $x_u(\cdot,x_0) = \{x_u(k;x_0),x_u(k+1;x_0),\dots,x_u(N;x_0) \}$ generated by the model~$x_u(k+1)=f_d(x_u(k),u(k))$ and stemming from~$x_0=x_u(0,x_0)$ satisfies~$x_u(k,x_0)\in \chi$ for all $k \in \{0,1,\dots,N\}$
    \end{enumerate}
\end{defn}
\end{comment}

\subsection{Dynamics Overview} \label{sec:dynamics}
\begin{figure}[ht]
\centering
\includegraphics[width=0.35\textwidth]{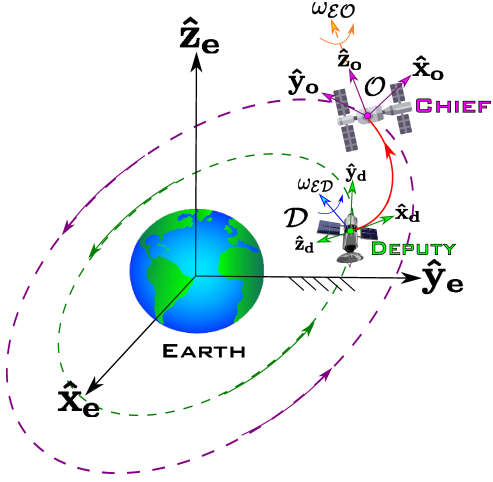}
\caption{The chief (with rotating orbit-fixed frame $\coordframe{O}$) and a deputy (with body-fixed frame $\coordframe{D}$) are orbiting about the Earth with inertial frame $\coordframe{E}$. The dashed lines are the closed orbital trajectories of both spacecraft. The red solid line depicts the rendezvous trajectory from the deputy to the chief as seen in the inertial frame.}
\label{fig:ExclusionZones}
\end{figure}

The objective of this work is to design an autonomous controller such that a deputy spacecraft's translational and attitudinal states converge to a successful docking configuration with the chief, as shown in Figure~\ref{fig:ExclusionZones}. 
Three frames of reference are used: 
\begin{enumerate}
\item {The Earth-Centered Inertial (ECI) frame originating from the Earth's center is denoted $\coordframe{E} := \{\mathbf{\hat{x}_e}, \mathbf{\hat{y}_e}, \mathbf{\hat{z}_e}\}$.}
\item {The Clohessy-Wiltshire (CW) frame is a non-inertial frame attached to the chief satellite, denoted $\mathcal{O} := \{\mathbf{\hat{x}_o}, \mathbf{\hat{y}_o}, \mathbf{\hat{z}_o}\}$. The origin lies at the center of mass of the chief, unless otherwise specified. The $\mathbf{\hat{x}_o}$-axis aligns with the vector pointing from the center of the Earth towards the center of mass of the chief satellite. The $\mathbf{\hat{z}_o}$-axis is in the direction of the orbital angular momentum vector of the chief spacecraft. The $\mathbf{\hat{y}_o}$-axis completes the right-handed orthogonal frame.}
\item{The deputy-fixed frame is a non-inertial body-fixed frame attached to the deputy satellite and denoted $\mathcal{D} := \{\mathbf{\hat{x}_d}, \mathbf{\hat{y}_d}, \mathbf{\hat{z}_d}\}$. The origin is assumed to lie at the center of mass of the deputy, unless otherwise specified. We use the convention that the $\mathbf{\hat{x}_d}$-axis is along the starboard direction, $\mathbf{\hat{y}_d}$ is in the fore direction, and $\mathbf{\hat{z}_d}$ completes the right-handed frame in the topside direction.}
\end{enumerate}
Let $(m_c, \inertia{C})$ and $(m_d, \inertia{D})$ be the mass-inertia pairs of the chief and deputy respectively. We make the following assumptions regarding their relative motion dynamics:

\begin{assumption} \label{as:A1}
    Both the chief and the deputy spacecraft are assumed to be rigid bodies of constant mass
\end{assumption}
\begin{assumption} \label{as:A2}
    The chief-fixed frame is assumed to always align with the orbit-fixed frame $\coordframe{O}$.
\end{assumption}
\begin{assumption} \label{as:A3}
    The chief spacecraft is in an uncontrolled, circular orbit.
\end{assumption}
\begin{assumption} \label{as:A4}
    The distance between the chief and the deputy is much less than the distance between the Earth and the chief.
\end{assumption}
\begin{assumption} \label{as:A5}
    The deputy spacecraft has bi-directional thrusters and external torque-generators installed along the axes aligning with its body-fixed frame.
\end{assumption}
\begin{assumption} \label{as:A6}
    The axes of the deputy-fixed frame are aligned with the principal axes of the deputy body.
\end{assumption}

\begin{comment}
    \begin{description}
\item[A1]{Both the chief and the deputy spacecraft are assumed to be rigid bodies of constant mass.}
\item[A2]{The chief-fixed frame is assumed to always align with the orbit-fixed frame $\coordframe{O}$.}
\item[A3]{The chief spacecraft is in an uncontrolled, circular orbit.}
\item[A4]{The distance between the chief and the deputy is much less than the distance between the Earth and the chief.}
\item[A5]{The deputy spacecraft has bi-directional thrusters and torque-generators installed along the axes aligning with its body-fixed frame.}
\item[A6]{The axes of the deputy-fixed frame are aligned with the principal axes of the deputy body.}
\end{description}
\end{comment}

 The position and velocity of the deputy relative to the chief with respect to the CW frame $\coordframe{O}$ are
\begin{equation}
\delta r^{\coordframe{O}} := (\delta x \ \delta y \ \delta z)^{\top} 
\textnormal{ and }
\delta \dot{r}^{\coordframe{O}} := (\delta \dot{x} \ \delta \dot{y} \ \delta \dot{z})^{\top}.
\end{equation}
Assumptions A3-A5 allow the use of the Clohessy-Wiltshire translational dynamics (see \citep{curtis2013orbital}) with deputy controls, expressed as
\begin{small}
\begin{equation}
\label{eq:CWdynamics_translational}
\begin{pmatrix}
\delta \dot{x} \\ \delta \dot{y} \\ \delta \dot{z} \\ \delta \ddot{x} \\ \delta \ddot{y} \\ \delta \ddot{z}
\end{pmatrix}
= \begin{pmatrix}
0 & 0 & 0 & 1 & 0 & 0 \\
0 & 0 & 0 & 0 & 1 & 0 \\
0 & 0 & 0 & 0 & 0 & 1 \\
3n^2& 0 & 0 & 0 & 2n & 0 \\
 0 & 0 & 0 & -2n & 0 & 0 \\
 0 & 0 & -n^2 & 0 & 0 & 0 
\end{pmatrix}
\begin{pmatrix}
\delta x \\ \delta y \\ \delta z \\ \delta \dot{x} \\ \delta \dot{y} \\ \delta \dot{z}
\end{pmatrix} + 
\begin{pmatrix}
\mathbf{0}_{3 \times 3}\\
\frac{1}{m_d}\rotmat{D}{O}
\end{pmatrix}
F^{\coordframe{D}}_d,
\end{equation}
\end{small}
 where the ``mean motion" constant is $n = \sqrt{\frac{\mu}{\norm{r_c}^3}}$ in $\text{rad}\cdot s^{-1}$, $\mu = 398600.4418 \ \text{km}^{3}s^{-2}$ is the Earth's standard gravitational parameter and $r_c$ is the radius of the chief's circular orbit.

\subsection{Attitude Dynamics - Chief Frame} \label{sec:new_rotation_dynamics}
Note the rotation $\rotmat{D}{O}$ applied to the deputy thrust vector $F^{\coordframe{D}}_d$ in \eqref{eq:CWdynamics_translational}. Both frames $\coordframe{O}$ and $\coordframe{D}$ are rotating with respect to the inertial frame $\coordframe{E}$. The angular velocity equations are
%\begin{align*}
$\angularvel{E}{D}{O} := (\omega_1 \ \omega_2 \ \omega_3)^{\top}$,  
$\angularvel{O}{D}{O} := (\delta\omega_1 \ \delta\omega_2 \ \delta\omega_3)^{\top}$,
and
$\angularvel{E}{O}{O} := (0 \ 0 \ n)^{\top}$. 
%\end{align*}
The deputy external torque vector in the $\coordframe{O}$ frame is $\tau^{\coordframe{O}} := (\tau_1 \ \tau_2 \ \tau_3)^{\top}$ and from assumption A6 the deputy inertia matrix in the deputy-fixed frame is $$\inertia{D} = \begin{pmatrix}
J_1 & 0 & 0 \\ 0 & J_2 & 0 \\ 0 & 0 & J_3
\end{pmatrix}.$$

The angular velocity of the $\coordframe{D}$ frame relative to the $\coordframe{O}$ frame is 
$\angularvel{O}{D}{E} = \angularvel{E}{D}{E} - \angularvel{E}{O}{E}$,
with the time derivatives given by 
\begin{equation}
\label{eq:rel_ang_vel_inertial}
\dotangularvel{O}{D}{E} = \dotangularvel{E}{D}{E} - \dotangularvel{E}{O}{E}.
\end{equation}
From \eqref{eq:rotation_velocity_a} and \eqref{eq:rotation_velocity_b} we yield the relations 
\begin{align*}
\dotangularvel{O}{D}{E} &= \rotmat{O}{E}\skewop{\angularvel{E}{O}{O}}\angularvel{O}{D}{O} + \rotmat{O}{E}\dotangularvel{O}{D}{O} \\
\dotangularvel{E}{O}{E} &= \rotmat{O}{E}\underbrace{\skewop{\angularvel{E}{O}{O}}\angularvel{E}{O}{O}}_{=0} + \rotmat{O}{E}\underbrace{\dotangularvel{E}{O}{O}}_{=0} = 0  \\
\dotangularvel{E}{D}{E} &= \rotmat{D}{E}\underbrace{\skewop{\angularvel{E}{D}{D}}\angularvel{E}{D}{D}}_{=0} + \rotmat{D}{E}\dotangularvel{E}{D}{D}.
\end{align*}
Substituting these into a manipulation of \eqref{eq:rel_ang_vel_inertial} gives
the relation 
%\begin{equation}
%\label{eq:relframe_quat}
$\dotangularvel{O}{D}{O} =  \rotmat{D}{O}\dotangularvel{E}{D}{D} - \skewop{\angularvel{E}{O}{O}}\angularvel{O}{D}{O}$.
%\end{equation}
%with accompanying relative quaternion kinematics
%\begin{align}
%\label{eq:relframe_eta}
%\dotquatscalar{D}{O} &= \frac{1}{2}[\quatvector{D}{O}]^{\top}\angularvel{O}{D}{O} \\
%\label{eq:relframe_rho}
%\dotquatvector{D}{O} &= -\frac{1}{2}[\quatscalar{D}{O}\identity{3} + \skewop{\quatvector{D}{O}}]\angularvel{O}{D}{O}. 
%\end{align}
Let the deputy body have inertia matrix $\inertia{D}$ as measured in the deputy-fixed frame with external control torque vector $\tau^{\coordframe{D}}$ being applied about the center of mass. From Euler's second law of motion we have
\begin{equation}
\label{eq:deputy_slew}    
\dotangularvel{E}{D}{D} = -\inertiainv{D}\skewop{\angularvel{E}{D}{D}}[\inertia{D}\angularvel{E}{D}{D}] + \inertiainv{D}\tau^{\coordframe{D}},
\end{equation}
where $\inertiainv{D} := [\inertia{D}]^{-1}$. Application of rotational transformation to \eqref{eq:deputy_slew} gives
$$\dotangularvel{O}{D}{O} =  \rotmat{D}{O}[-\inertiainv{D}\skewop{\angularvel{E}{D}{D}}[\inertia{D}\angularvel{E}{D}{D}] + \inertiainv{D}\tau^{\coordframe{D}}] - \skewop{\angularvel{E}{O}{O}}\angularvel{O}{D}{O},$$
which expands to
\begin{dmath}
\label{eq:ugly_ass_equation}
\dotangularvel{O}{D}{O} = \skewop{\angularvel{O}{D}{O}}\angularvel{E}{O}{O} + \inertiainv{O}(t)\tau^{\coordframe{O}} -  \inertiainv{O}(t)\big(\skewop{\angularvel{O}{D}{O}}[\inertia{D}\rotmat{O}{D}\angularvel{O}{D}{O}] + \skewop{\angularvel{O}{D}{O}}[\inertia{D}\rotmat{O}{D}\angularvel{E}{O}{O}] + \skewop{\angularvel{E}{O}{O}}[\inertia{D}\rotmat{O}{D}\angularvel{O}{D}{O}] + \skewop{\angularvel{E}{O}{O}}[\inertia{D}\rotmat{O}{D}\angularvel{E}{O}{O}]\big).
\end{dmath}
Here $\inertiainv{O} = \rotmat{D}{O}\inertiainv{D}\rotmat{O}{D}$ is a time-varying inertia matrix.

\subsection{Docking Configuration} \label{sec:docking_config}
To achieve a successful docking configuration, the deputy must reach a set of predefined relative translational and attitudinal states. Let $\quatdes{D}{O}$ be the quaternion representing the desired rotation from the deputy $\coordframe{D}$ frame to the CW $\coordframe{O}$ frame, and let $\quat{D}{O}$ be the actual (i.e., plant-state) rotation quaternion. Analogously, let $\angularveldes{O}{D}{O}$ be the desired relative angular velocity and $\angularvel{O}{D}{O}$ be the actual relative angular velocity. Lastly, let $\quaterr{D}{O} := \{\quatdes{D}{O}\}^{-1}\otimes\quat{D}{O}$ and $\angularvelerr{O}{D}{O} := \angularvel{O}{D}{O} - \angularveldes{O}{D}{O}$ respectively denote the error quaternion and error angular velocity. 

To achieve on-orbit attitudinal synchronization between both spacecraft, the aim is to drive $\quaterr{D}{O} \to q^{I}$ and $\angularvelerr{O}{D}{O} \to \mathbf{0}_3$. In the ARPOD problem, $\quatdes{D}{O} = q^{I}$ and $\angularveldes{O}{D}{O} = \mathbf{0}_3$, which results in the attitudinal error 
kinematics~$\dotquatscalarerr{D}{O} = \frac{1}{2}[\quatvectorerr{D}{O}]^{\top}\angularvelerr{O}{D}{O}$
and~$\dotquatvectorerr{D}{O} = -\frac{1}{2}[\quatscalarerr{D}{O}\identity{3} + \skewop{\quatvectorerr{D}{O}}]\angularvelerr{O}{D}{O}$. 
%
%\begin{align}
%\label{eq:error_eta}
%\dotquatscalarerr{D}{O} &= \frac{1}{2}[\quatvectorerr{D}{O}]^{\top}\angularvelerr{O}{D}{O} \\
%\label{eq:error_rho}
%\dotquatvectorerr{D}{O} &= -\frac{1}{2}[\quatscalarerr{D}{O}\identity{3} + \skewop{\quatvectorerr{D}{O}}]\angularvelerr{O}{D}{O}. 
%\end{align}
Relative translational states required for docking are $\delta r^{\coordframe{O}} = \textbf{0}_{3}$ and $\delta \dot{r}^{\coordframe{O}} = \textbf{0}_{3}$. The state vector of interest is
\begin{equation}
x \coloneqq \Big([\delta r^{\coordframe{O}}]^{\top}, \ [\delta\dot{r}^{\coordframe{O}}]^{\top}, \ \quatscalarerr{D}{O}, \ [\quatvectorerr{D}{O}]^{\top}, \ [\angularvelerr{O}{D}{O}]^{\top}\Big)^{\top} \in \chi
\end{equation}
and the objective docking state be
\begin{equation}
\label{eq:docking}
x_d \coloneqq  ([\mathbf{0}_3]^{\top}, \ [\mathbf{0}_3]^{\top}, \ 1, \ [\mathbf{0}_3]^{\top}, \ [\mathbf{0}_3]^{\top})^{\top} \in \chi.
\end{equation}
This vector~$x_d$ 
indicates that the frame~$\coordframe{O}$ has the same origin, linear velocity, orientation, and angular velocity as frame~$\coordframe{D}$. 
Let the control vector of interest in frame $\coordframe{O}$ be
%\begin{equation}
%\label{eq:inputs}
$u \coloneqq  ([F^{\coordframe{D}}_d]^{\top}, \ [\tau^{\coordframe{D}}]^{\top})^{\top} \in \mathcal{U}$
%\end{equation}
and the objective inputs be
%\begin{equation}
%\label{eq:inputsDesired}
$u_d \coloneqq  ([\mathbf{0}_3]^{\top}, \ [\mathbf{0}_3]^{\top})^{\top} \in \mathcal{U}$.
%\end{equation}
Let us define~$z(k) \coloneqq [u(k)^\top \ x(k)^\top]^\top$ and $z_d \coloneqq [u_d^\top \ x_d^\top]^\top$. The next section uses this model to solve the ARPOD problem using
time-constrained MPC.

%===================================================================================================================

\section{Main Results} \label{sec:results}
In this section we present two simulations. First, we consider time-constrained MPC with unperturbed states. Second, we consider states affected by perturbations to illustrate the robustness of time-constrained MPC.

\subsection{Time-Constrained MPC for ARPOD}

The problem we solve in this section is given next. 

\begin{problem}[Time-Constrained MPC for ARPOD] \label{prob:arpodMPC}
%\mh{Shouldn't the~$x_d$ in the cost have some dependence on~$k$?}
%\gb{The docking configuration~$x_d$ is constant the entire problem horizon as in Equation~\eqref{eq:docking}}
\begin{small}
\begin{equation} 
    \begin{split}
        \underset{\textbf{u}(k)}{\textnormal{minimize}} &\ J(x(k),\textbf{u}(k)) = \\ &\ \sum^{N-1}_{i=k} (x(i)-x_d)^\top Q (x(i)-x_d) + u(i)^\top R u(i) \\
        \textnormal{subject to} &\ x(k) = x(0) \\
                                &\ x(N) = x_d \\
                                &\  x(i+1) = g_d(x(i),u(i)), \quad i=k,k+1,\dots,N-1 \\
                                &\ u_\min \leq u(i) \leq u_\max, \quad i=k,k+1,\dots,N-1\\
                                &\ \textbf{u}(k)\in \U_N \\
                                & \ j_k \leq j_\max,
    \end{split}
\end{equation}
\end{small}
where~$Q \in \R^{13\times 13}$, ~$R \in \R^{6\times 6}$,~$u_\min \coloneqq (-10^{-3} \times \ones^\top_3, \ -10^{-4} \times \ones^\top_3)^\top$,~$u_\max \coloneqq -u_\min$, and~$g_d:\chi \times \U \rightarrow \chi$ are the discretized versions 
of~\eqref{eq:quaternion_velocity_c},\eqref{eq:CWdynamics_translational}, and \eqref{eq:ugly_ass_equation} generated by the fourth-order Runge-Kutta method.
\end{problem}

%The continuous-time dynamics were discretized using the fourth order Runge-Kutta method.

\subsection{Unperturbed Results}
We consider the problem parameters and initial states in Table~\ref{tb:parameters}, where~$n$ is the chief mean motion, $m_d$ is the mass of the deputy, and~$J_d$ is its moment of inertia matrix. The values in Table~\ref{tb:parameters} represent a chief spacecraft in low Earth orbit and a deputy in 
a typical ARPOD initial condition.

\begin{table}[htb]
\caption[]{\label{tb:parameters} Problem Parameters and Initial State}
\begin{center}
\begin{tabular}{ |c|c|c|c| } 
\hline
Parameter & Value & Units \\
\hline
$n$ & -0.0011 & rad/s \\
$m_d$ & 12 & kg \\ 
$J_d$ & diag([0.2734, 0.2734, 0.3125]) & kg$\cdot \text{m}^2$ \\ 
$\angularvel{E}{O}{O}$ & $(0 \ 0 \ -0.0011)^T$ & rad/s \\
$\delta r^{\coordframe{O}}$ & $(1.5 \ -1.77 \ 3)^T$ & km \\
$\delta \dot{r}^{\coordframe{O}}$ & $(0.001 \ 0.0034 \ 0)^T$ & km/s \\
$q^{\coordframe{O}}_{\coordframe{D}}$ & $(0.7715 \ 0.4629 \ 0.3086 \ 0.3086)^T$ &  \\
$\omega^{\coordframe{O}}_{\coordframe{O}\coordframe{D}}$ & $(0 \ 0 \ -0.005)^T$ & rad/s \\
\hline
\end{tabular}
\end{center}
\end{table}
 The cost matrices used were~$Q=\textnormal{diag}([10\times \ones^\top_3,  10^{-4}\times \ones^\top_3, 10^8 \times \ones^\top_7])$ and $R=\textnormal{diag}([10^3\times \ones^\top_3, 10^{10} \times \ones^\top_3])$. 
 We considered the prediction horizon~${N=1000}$ and sampling time~${t_s = 3 s}$. 
 Since the torque inputs and attitudinal states are heavily penalized, the orientation of the deputy is prioritized over the translational state. 
 We formulated Problem~\ref{prob:arpodMPC} 
 in MATLAB using the CasADi symbolic framework and solved it using Ipopt~\citep{wachter2006}. 
 
 In the minimization of Problem~\ref{prob:arpodMPC}, at each time~$k$ we designed two stopping conditions, either of which terminates 
 the minimization if it is met: 
 (1) the change in 
 objective function value at consecutive iterations was less than a specified tolerance or (2) the algorithm had completed the maximum allowable number of iterations, i.e.,~$j_k = j_\max$, to enforce the computational time constraint. In this subsection, we declare that 
 we have achieved the docking configuration when~${\Vert z(k) - z_d \Vert_\infty \leq 10^{-3}}$. 
 %\mh{Have we defined~$z(k)$ and~$z_d$ here? Make sure they're defined explicitly.}
 %\gb{I defined~$z(k)$ at the end of the docking configuration section to be $z(k) \coloneqq [u(k)^\top \ x(k)^\top]^\top$ and $z_d \coloneqq [u_d^\top \ x_d^\top]^\top$. Maybe there is a better place for them.}
 
\begin{figure}
    \centering
    \includegraphics[width = 0.40 \textwidth]{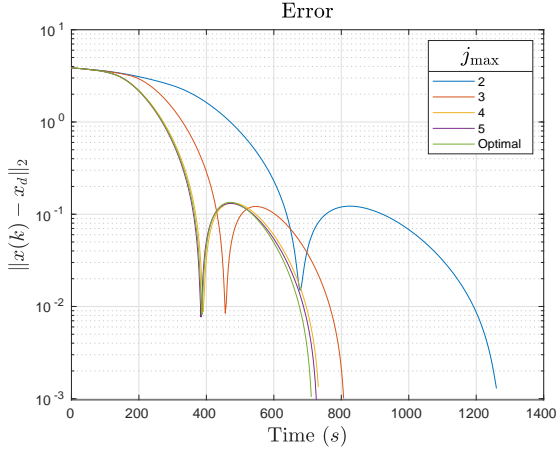}
    \caption{Plot of the unperturbed error~$\Vert x(k) - x_d \Vert_2$ 
    as a function of time for various values of~$j_{max}$. 
    %is shown for various values of the maximum
    %allowable number of iterations,~$j_{max}$. For~$j_{max} \geq 5$,
    %the resulting errors over time
    %are nearly indistinguishable from the optimal errors. 
    }
    \label{fig:unperturbError}
\end{figure}

 Figure \ref{fig:unperturbError} shows the distance between the current state~$x(k)$ and docking configuration~$x_d$ over time for multiple values of~$j_\max$. 
 The error~${\Vert x(k) - x_d \Vert_2}$ decreases faster as we allow more iterations to occur, which agrees with the intuition that more computations drive
 an algorithm's iterates closer to an optimum. 
 We can observe that the optimal trajectory, i.e., the one computed with~$j_\max = \infty$, achieves the docking state in the shortest amount of time.
 
\begin{figure}
    \centering
    \includegraphics[width = 0.40 \textwidth]{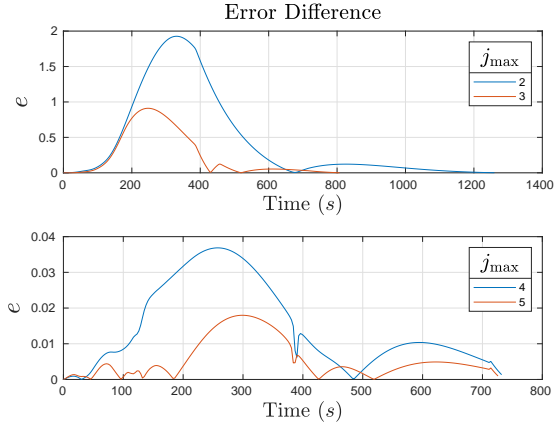}
    \caption{Comparison plot of the difference between time-constrained and optimal unperturbed error as a function of time. The value of~${e \coloneqq \vert \Vert \tilde{x}(k) - x_d \Vert_2 - \Vert x^*(k) - x_d \Vert_2 \vert}$ is shown for 
    various values of~$j_{max}$. 
    %the maximum number of allowed iterations,~$j_{max}$.
    %We see that the error incurred shrinks as~$j_{max}$ grows
    %from~$2$ to~$5$, and for~$j_{max} = 5$ the error always remains
    %below~$0.02$. 
    }
    \label{fig:unperturbErrorDiff}
\end{figure}

\begin{figure}
        \centering
        \begin{subfigure}[b]{0.23\textwidth}
            \centering
            \includegraphics[width=\textwidth]{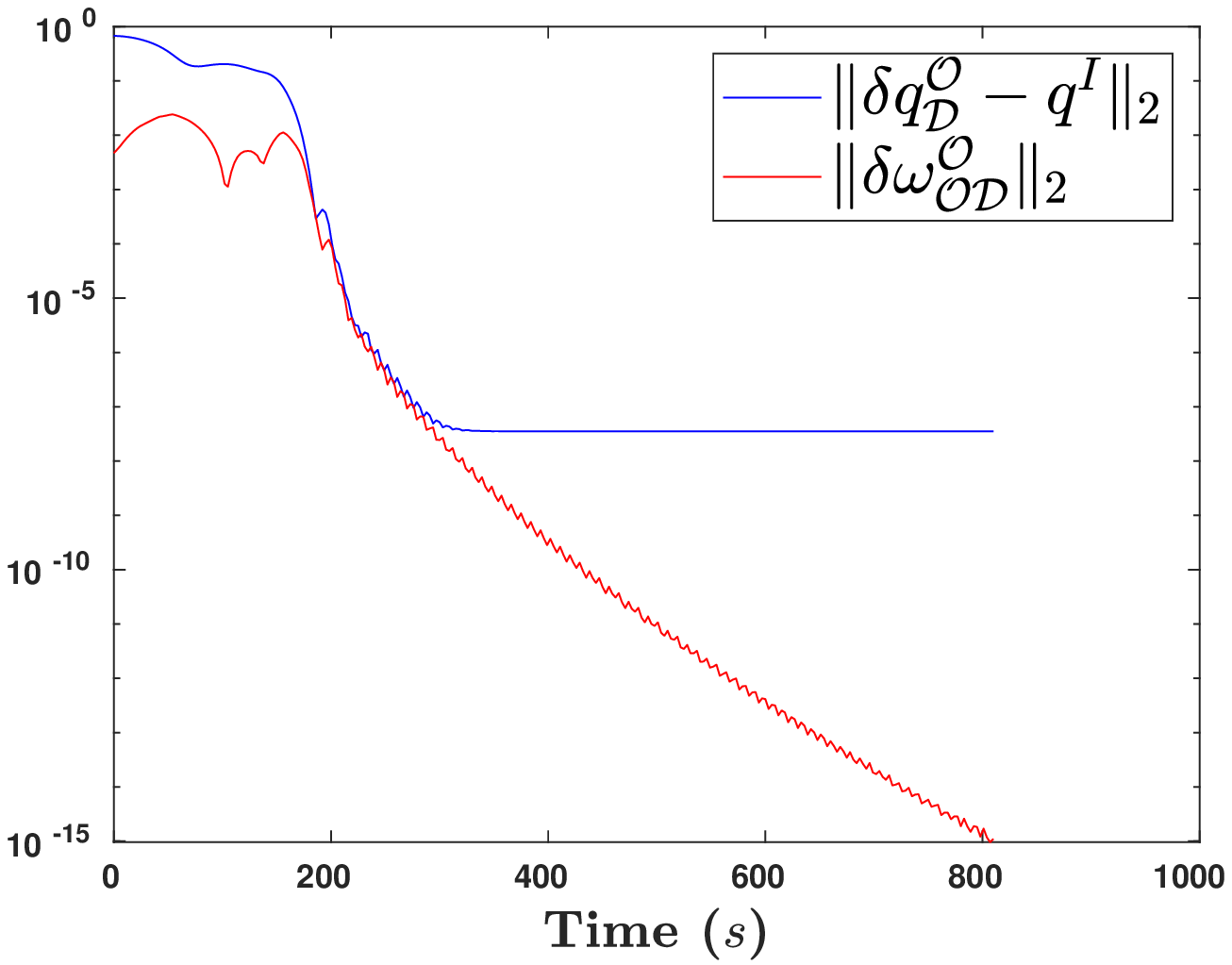}
            \caption[]%
            {{Attitude Error}}    
            \label{fig:attitudeError}
        \end{subfigure}
        \hfill
        \begin{subfigure}[b]{0.23\textwidth}  
            \centering 
            \includegraphics[width=\textwidth]{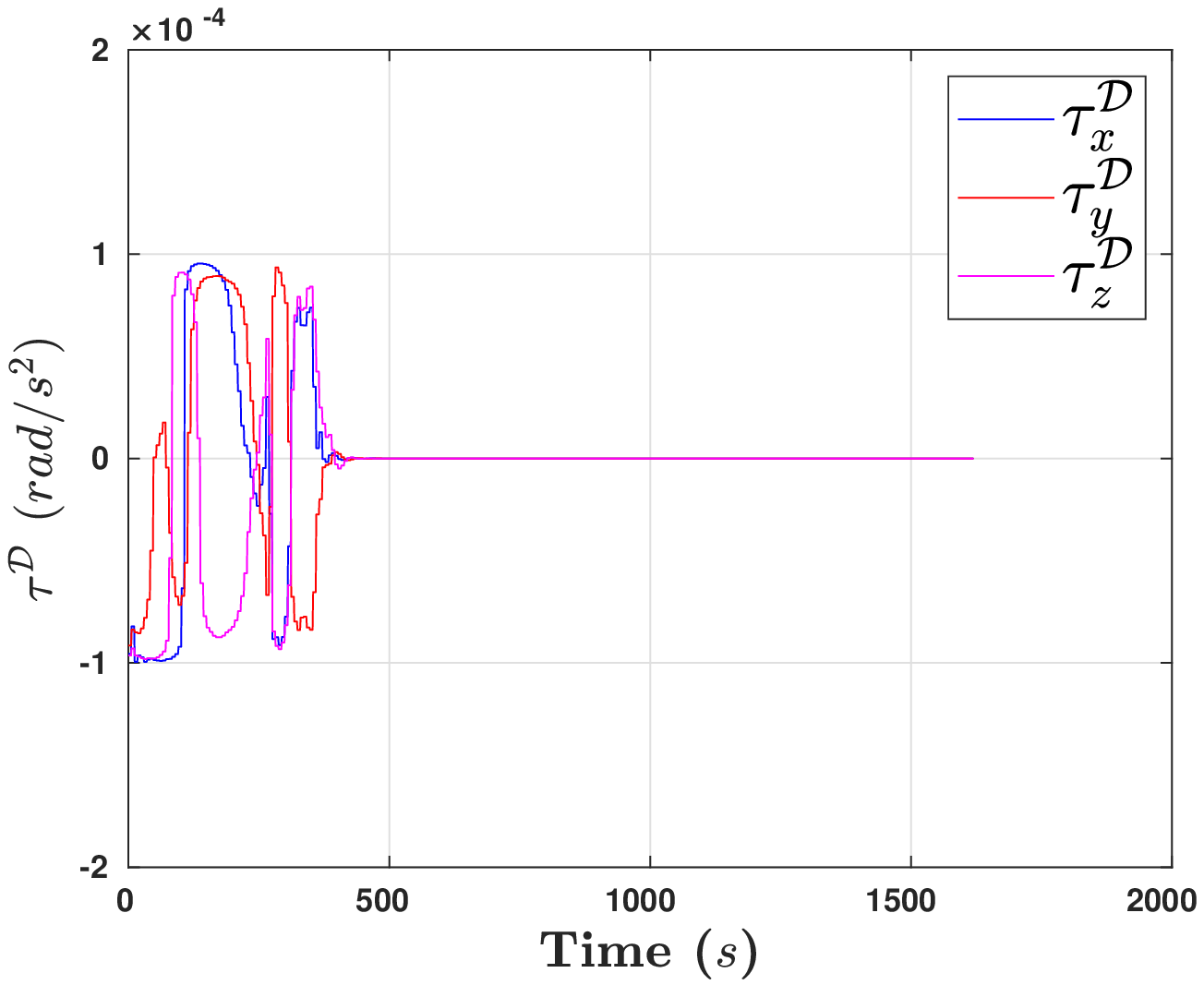}
            \caption[]%
            {{Torque Inputs,~$\tau^{\coordframe{D}} (\text{rad/s}^2$}) }    
            \label{fig:torque}
        \end{subfigure}
        \vskip\baselineskip
        \begin{subfigure}[b]{0.23\textwidth}   
            \centering 
            \includegraphics[width=\textwidth]{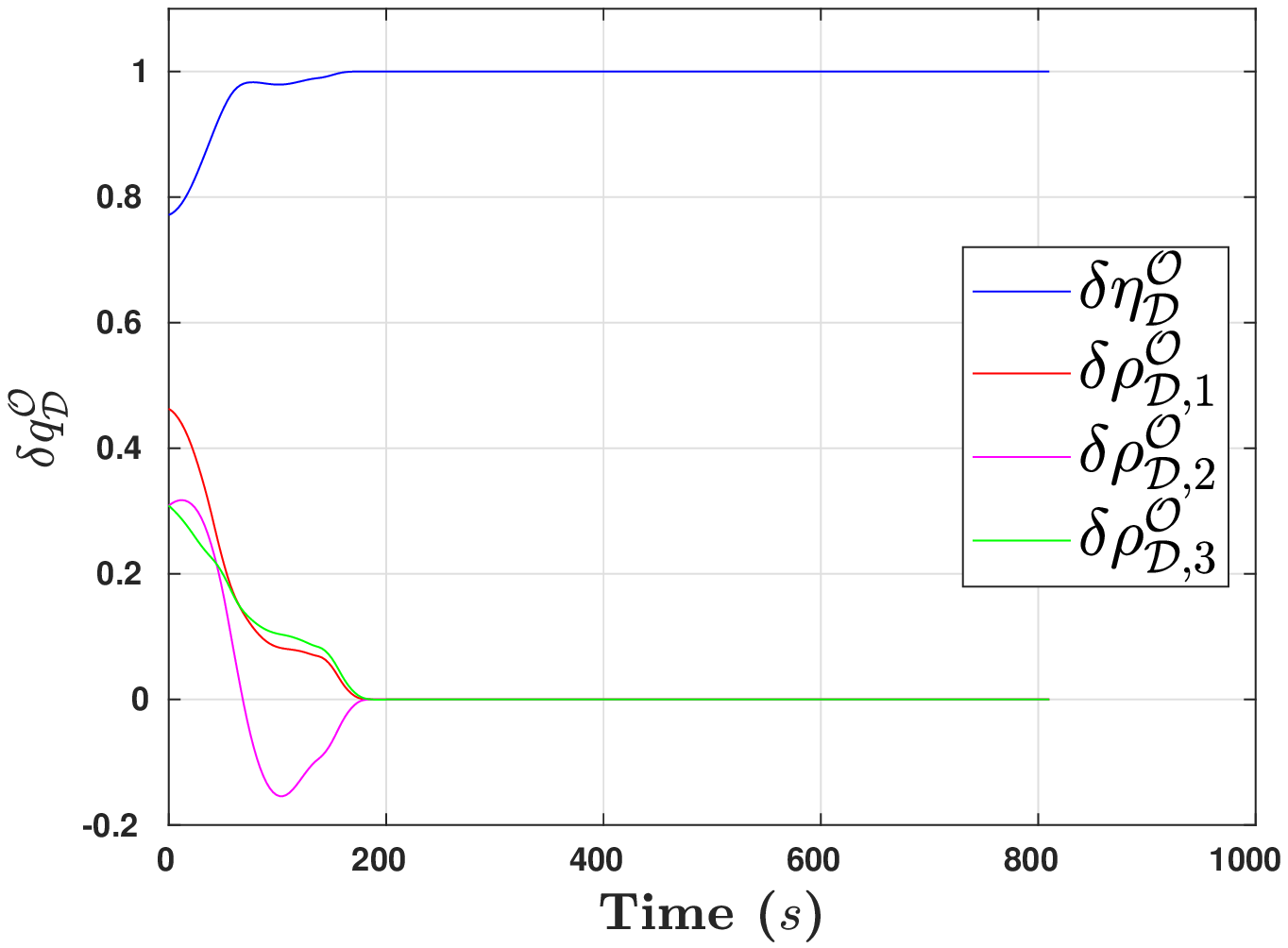}
            \caption[]%
            {{Error Quaternion,~$\delta q^{\coordframe{O}}_{\coordframe{D}}$}}    
            \label{fig:quaternion}
        \end{subfigure}
        \hfill
        \begin{subfigure}[b]{0.24\textwidth}   
            \centering 
            \includegraphics[width=\textwidth]{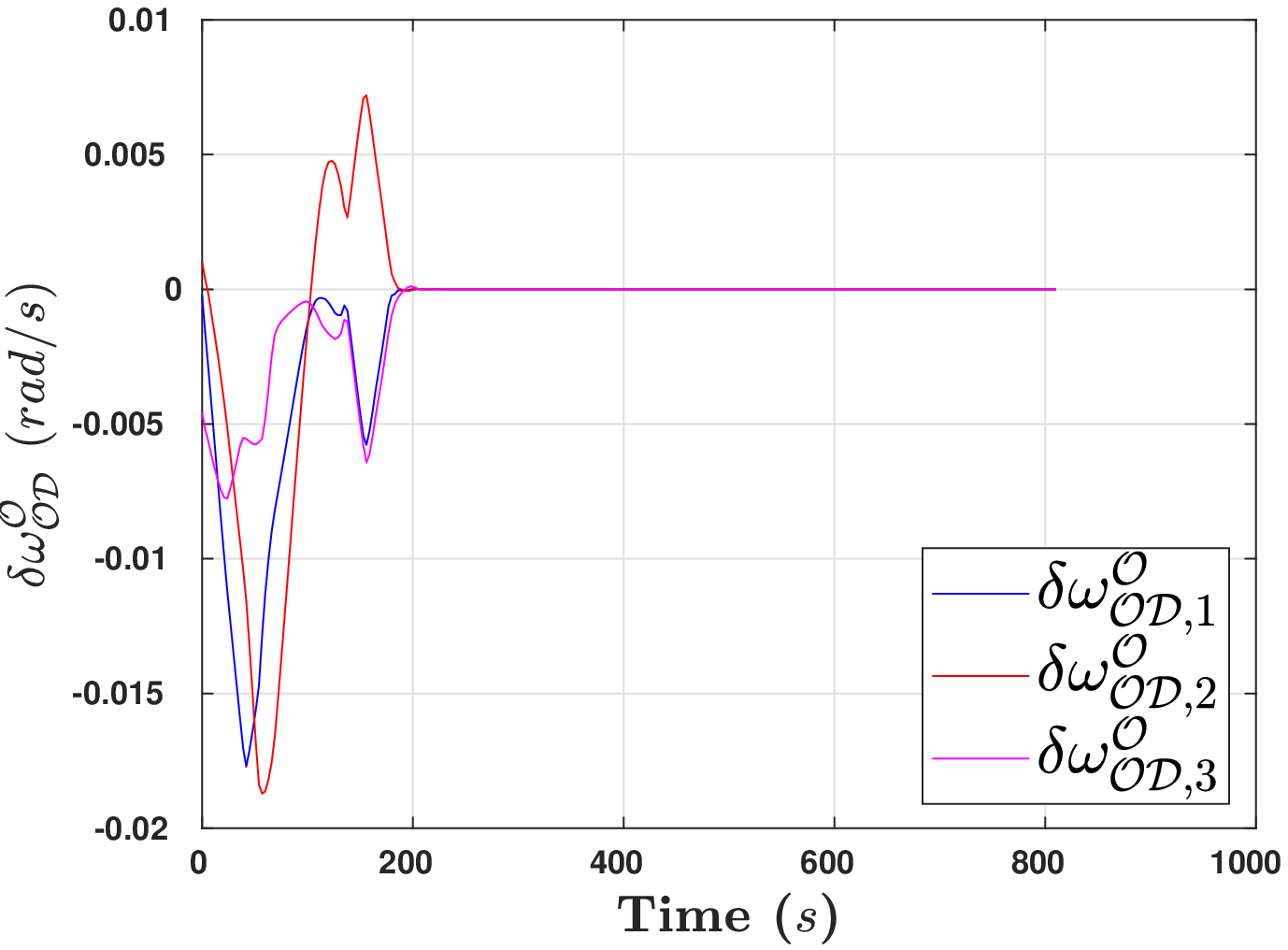}
            \caption[]%
            {{Error Angular Velocity,~$\delta \omega^{\coordframe{O}}_{\coordframe{O}\coordframe{D}}$ }}    
            \label{fig:angularVel}
        \end{subfigure}
        \caption[]
        {Plots of the unperturbed attitudinal error, states, and torque inputs 
        for~$j_{max} = 3$.} 
        \label{fig:attitude}
\end{figure}

\begin{figure}
        \centering
        \begin{subfigure}[b]{0.23\textwidth}
            \centering
            \includegraphics[width=\textwidth]{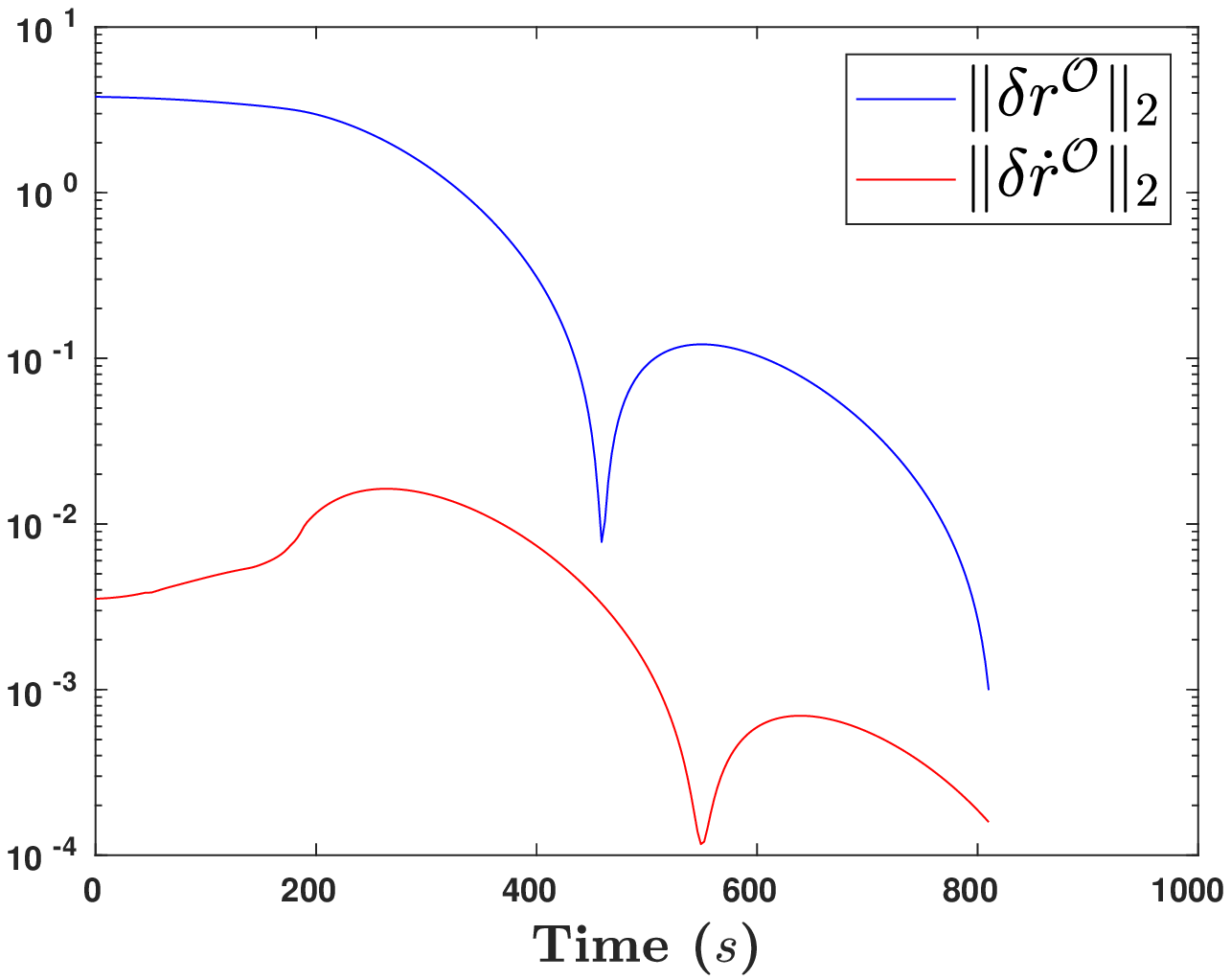}
            \caption[]%
            {{Translational Error}}    
            \label{fig:translationError}
        \end{subfigure}
        \hfill
        \begin{subfigure}[b]{0.23\textwidth}  
            \centering 
            \includegraphics[width=\textwidth]{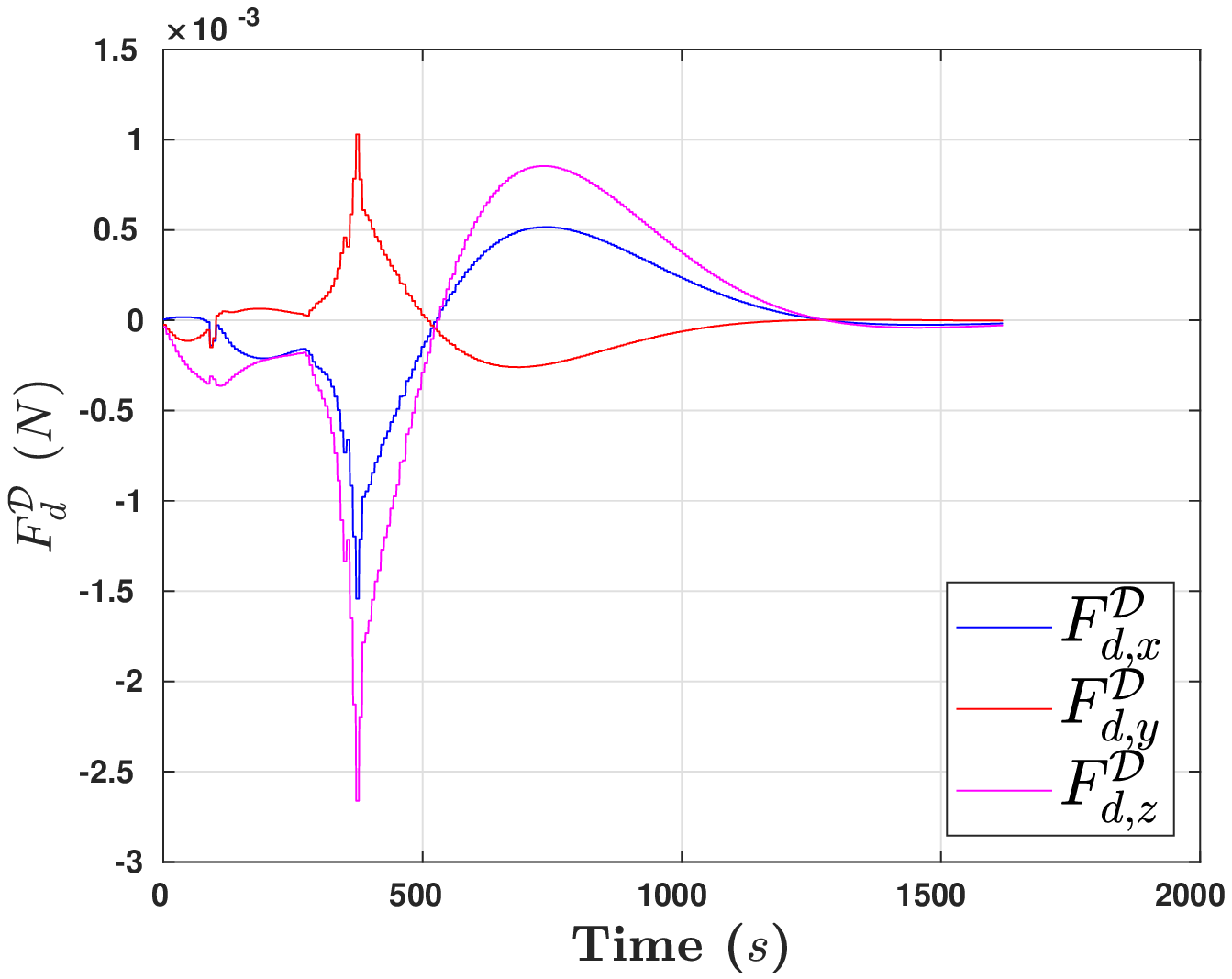}
            \caption[]%
            {{Thrust Inputs,~$F^{\coordframe{D}}_d$ (N)}}    
            \label{fig:thrust}
        \end{subfigure}
        \vskip\baselineskip
        \begin{subfigure}[b]{0.23\textwidth}   
            \centering 
            \includegraphics[width=\textwidth]{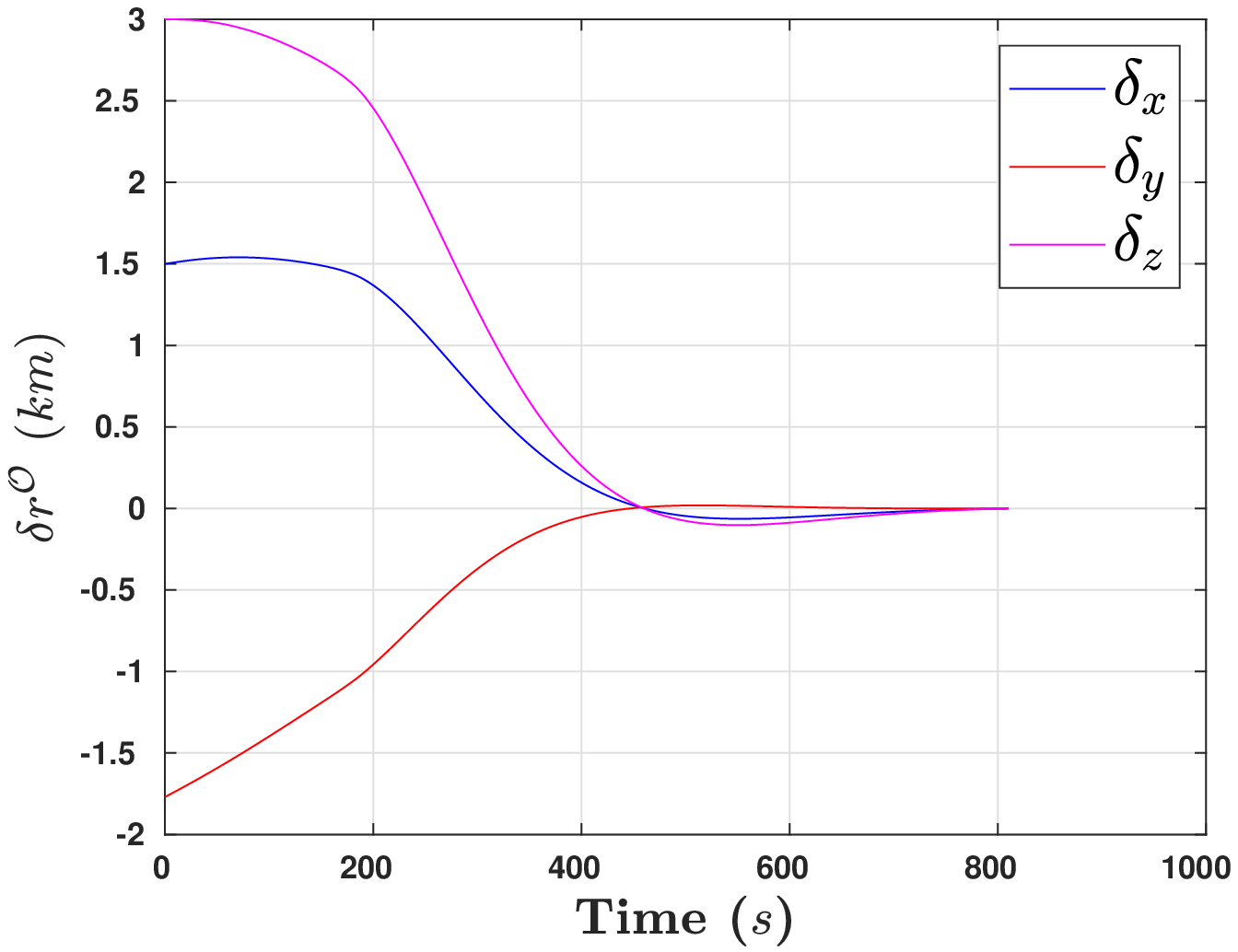}
            \caption[]%
            {{Position,~$\delta r^{\coordframe{O}}$ (km)}}    
            \label{fig:position}
        \end{subfigure}
        \hfill
        \begin{subfigure}[b]{0.23\textwidth}   
            \centering 
            \includegraphics[width=\textwidth]{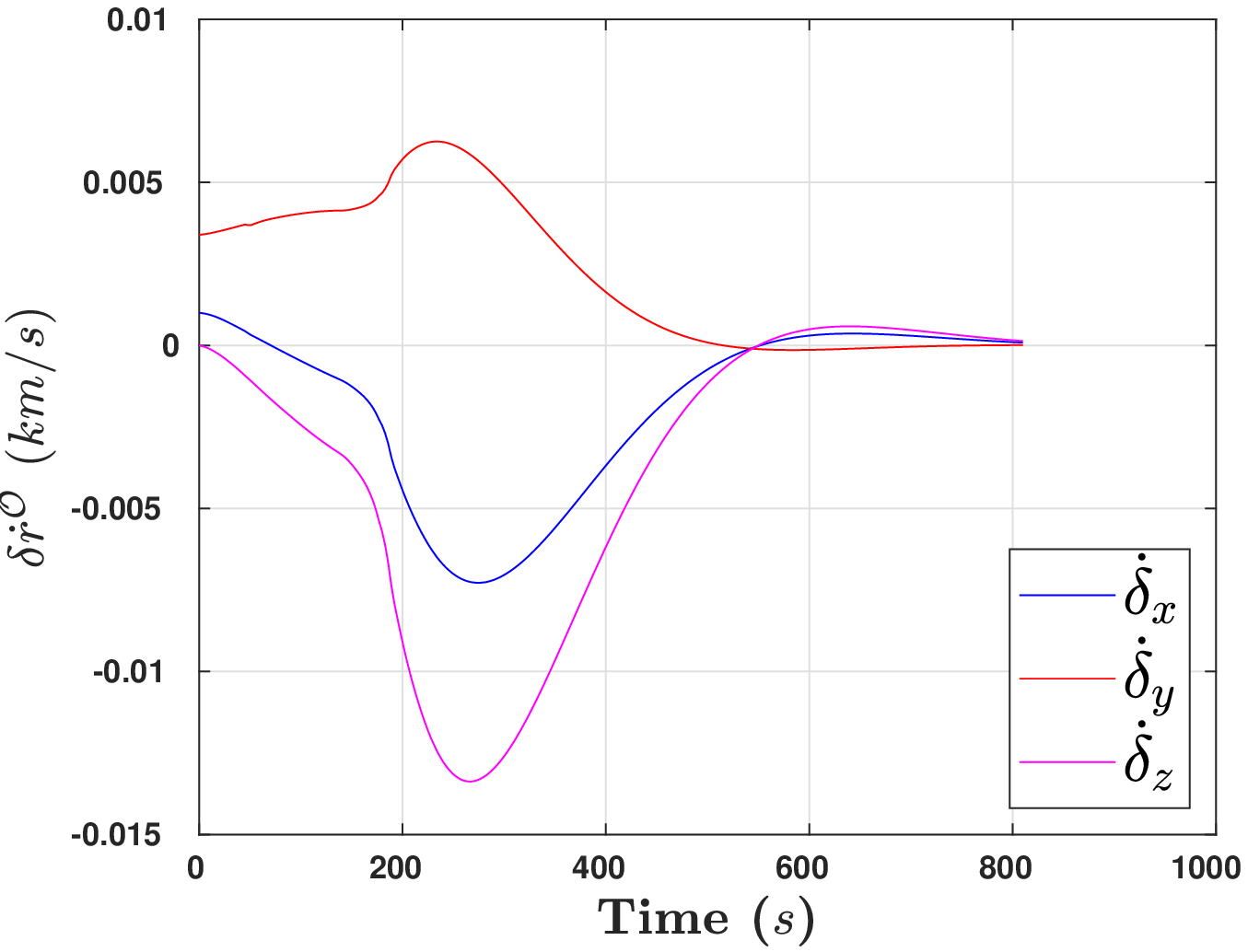}
            \caption[]%
            {{Linear Velocity,~$\delta \dot{r}^{\coordframe{O}}$ (km/s)}}    
            \label{fig:velocity}
        \end{subfigure}
        \caption[]
        {Plots of the unperturbed translational error, states, and thrust 
        inputs for~$j_{max}=3$.} 
        \label{fig:translation}
\end{figure}

\begin{figure}
    \centering
    \includegraphics[width = 0.40 \textwidth]{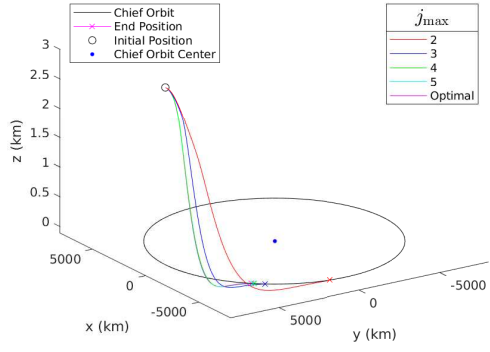}
    \caption{
    Plot of the the position of the deputy docking with the target is 
    shown for the unperturbed case for various 
    values of~$j_{max}$.
    }
    \label{fig:ECIorbits}
\end{figure}

 Furthermore, there is a trade-off between the value of~$j_\max$ and time it takes to reach~$x_d$.
 This is shown in Figure~\ref{fig:unperturbErrorDiff}, where we plot the absolute difference between the time-constrained and optimal error.
 In addition, Figures~\ref{fig:attitude} and~\ref{fig:translation} show 
 that 
 the attitudinal and translational states achieve the docking configuration given the initial conditions in Table~\ref{tb:parameters} with~$j_\max = 3$.
 Figure~\ref{fig:ECIorbits} shows the positional trajectories in 3-D space achieving the docking configuration.
 A Monte Carlo simulation was performed with~$j_\max = 3$ and the average error 
 for~$50$ initial conditions was $2.6\cdot 10^{-3}$. 
 These results show that time-constrained MPC can control the deputy to the docking configuration using unperturbed states. 
 %\mh{Gabe, for comparison, can we say how many iterations are needed  if we don't constrain~$j_k$? This is a key piece of the story that we want  to tell.  }
 %\gb{The number of iterations needed changes every $k$. Initially it is larger because we start at a random initial point. However, after about 5-10 MPC loops the number of iterations needed for termination are very similar for all cases of~$j_\max$ because the initial point used in the next MPC loop is "warm-started" by using the solution of the previous loop}
 Table~\ref{tb:unperturbTimes} shows that we achieve at least a 23.98\% reduction in average time to solve Problem~\ref{prob:arpodMPC} at each time~$k$, for all values of~$j_\max$ that we used. Moreover, Table~\ref{tb:unperturbTimes} shows a 
 reduction of more than 92\% in the maximum loop time
 for the worst case solving time for all~$j_\max$.
  Next, we demonstrate the robustness of time-constrained MPC to perturbations.
 
\begin{table} 
\caption[]{Solver times and reductions for Problem~\ref{prob:arpodMPC} without perturbations}
\resizebox{\columnwidth}{10mm}{
\begin{tabular}{|| M{1cm} | M{2.2cm}| M{2.3cm} | M{2.2cm} | M{2.3cm} ||} 
 \hline
 $j_{\max}$  &   Average Time/Loop Reduction (s)&  Average Time/Loop Reduction (\%) & Maximum Loop Time Reduction (s) & Maximum Loop Time Reduction (\%)\\  [0.5ex] 
 \hline\hline
 2 & 0.5159 & 30.52\% & 26.40 & 95.59\% \\
 \hline
 3 & 0.4245  & 25.11\% & 26.09 & 94.48\%\\
 \hline
 4 & 0.4283 & 25.33\% & 25.75 & 93.26\% \\
 \hline
 5 & 0.4054 & 23.98\% & 25.43 & 92.09\% \\ [1ex] 
 \hline
\end{tabular}
}
\label{tb:unperturbTimes}
\end{table}

\subsection{Perturbed Results}
\begin{figure}
    \centering
    \includegraphics[width = 0.40 \textwidth]{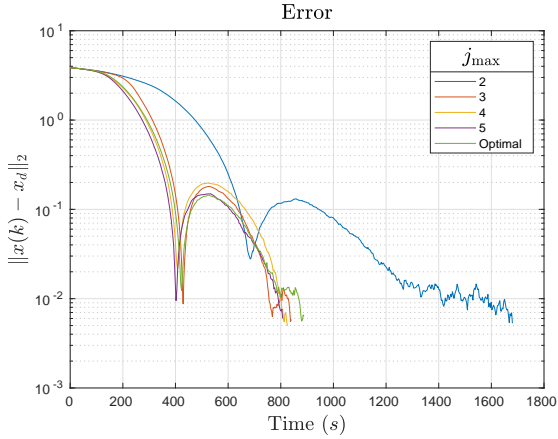}
    \caption{Plot of the perturbed error,~$\Vert x(k) - x_d \Vert_2$,
    as a function of time for various values of~$j_{max}$. 
    }
    \label{fig:perturbError}
\end{figure}

One advantage of MPC is robustness to perturbations. 
In these simulations, at each time~$k$, we add a disturbance~$\omega (k)$ to the state~$x(k)$. The initial state for the next MPC loop, at time~$k+1$ was
%\begin{equation}
    $x(0) = g_d(x(k),u(k)) + \omega(k)$.
%\end{equation}
The components of the perturbation vector ~$\omega(k) \in \R^{13}$ were sampled from  multivariate normal distributions as
\begin{equation}
    \omega(k) \sim \begin{pmatrix}
        10^{-3}\times \multiGauss{\mathbf{0}_{3}}{\identity{3}} \\
        10^{-6}\times \multiGauss{\mathbf{0}_{3}}{\identity{3}} \\
        10^{-8}\times \multiGauss{\mathbf{0}_{4}}{\identity{4}} \\
        10^{-6}\times \multiGauss{\mathbf{0}_{3}}{\identity{3}} \\
    \end{pmatrix}. 
\end{equation}
 We consider the same initial condition from Table~\ref{tb:parameters} and same problem parameters from the previous subsection. In this subsection, we declare that we have achieved the docking configuration when~$\Vert z(k) - z_d \Vert_\infty \leq 5\times 10^{-3}$. 
 %\mh{Same as above: have we defined~$z(k)$ and~$z_d$ here?}
 %\gb{I defined~$z(k)$ at the end of the docking configuration section to be $z(k) \coloneqq [u(k)^\top \ x(k)^\top]^\top$ and $z_d \coloneqq [u_d^\top \ x_d^\top]^\top$. Maybe there is a better place for them.}
 Figure~\ref{fig:perturbError} shows that the deputy still achieves the docking configuration with the chief subject to perturbations and limited computational time.
It takes more time to achieve the docking configuration for each~$j_\max$ compared to the unperturbed case (cf. Figure~\ref{fig:unperturbError}). 

\begin{figure}
        \centering
        \begin{subfigure}[b]{0.23\textwidth}
            \centering
            \includegraphics[width=\textwidth]{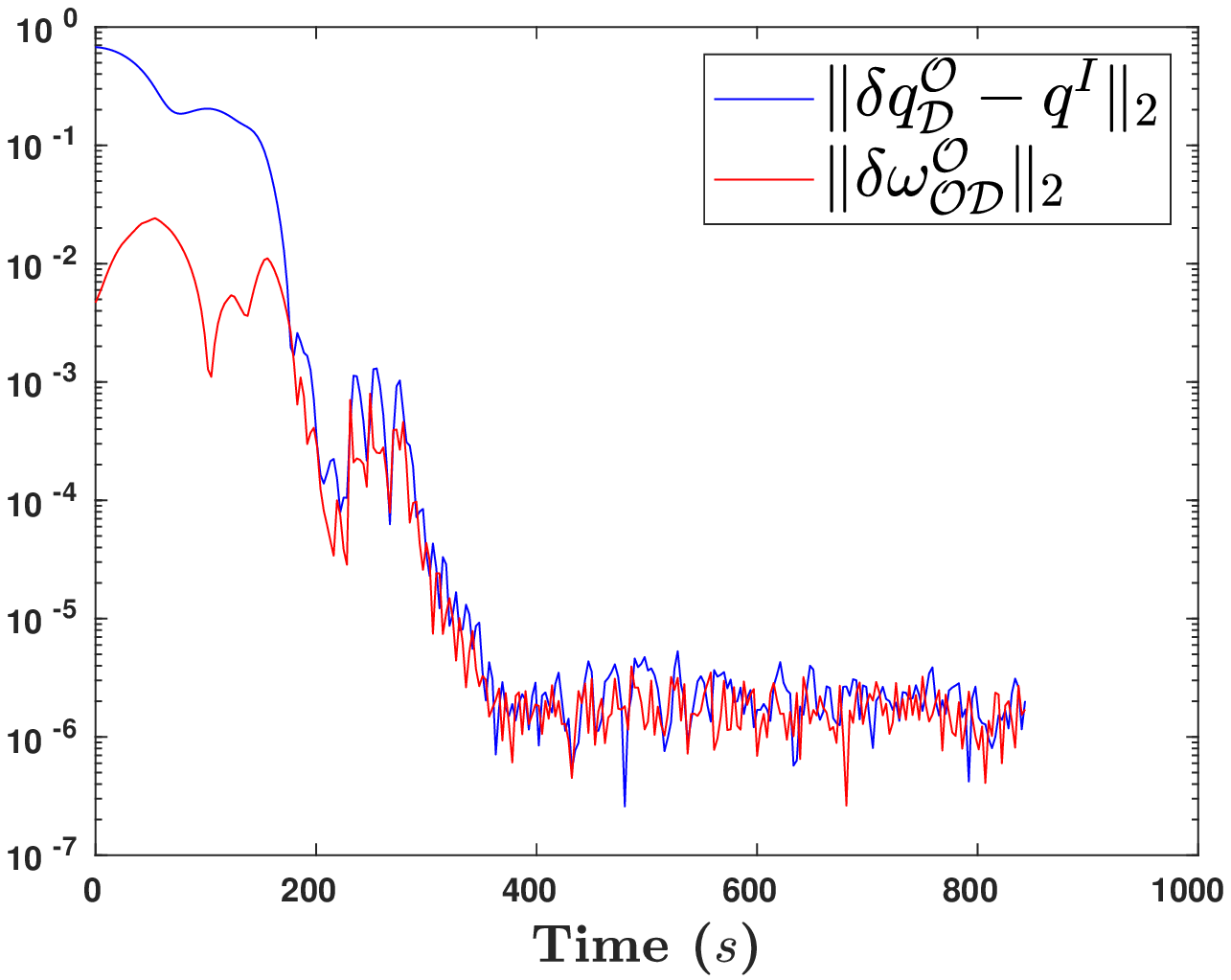}
            \caption[]%
            {{\small Attitude Error}}    
            \label{fig:attitudeError2}
        \end{subfigure}
        \hfill
        \begin{subfigure}[b]{0.23\textwidth}  
            \centering 
            \includegraphics[width=\textwidth]{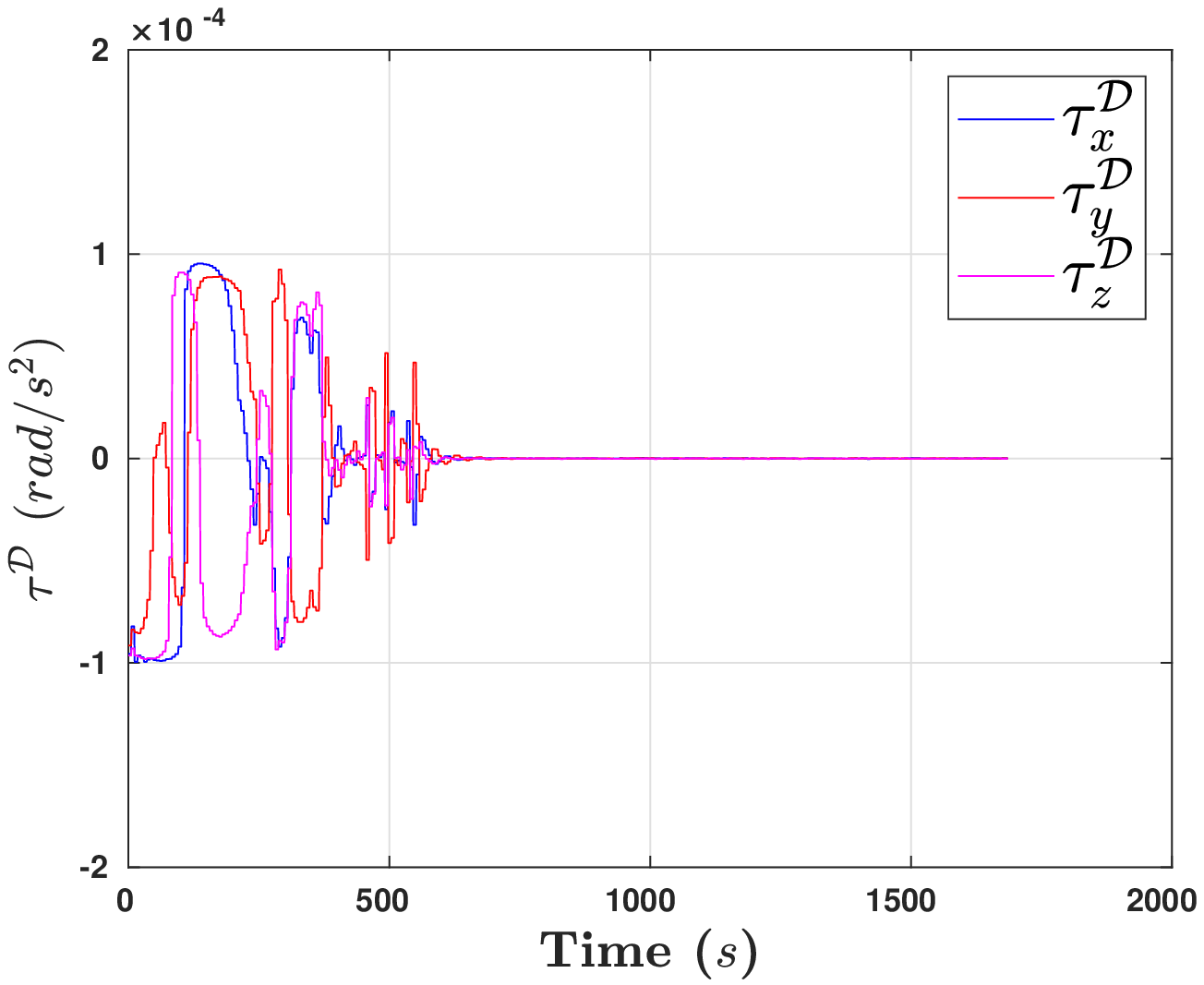}
            \caption[]%
            {{\small Torque Inputs,~$\tau^{\coordframe{D}}$ ($\text{rad/s}^2$)}}    
            \label{fig:torque2}
        \end{subfigure}
        \vskip\baselineskip
        \begin{subfigure}[b]{0.23\textwidth}   
            \centering 
            \includegraphics[width=\textwidth]{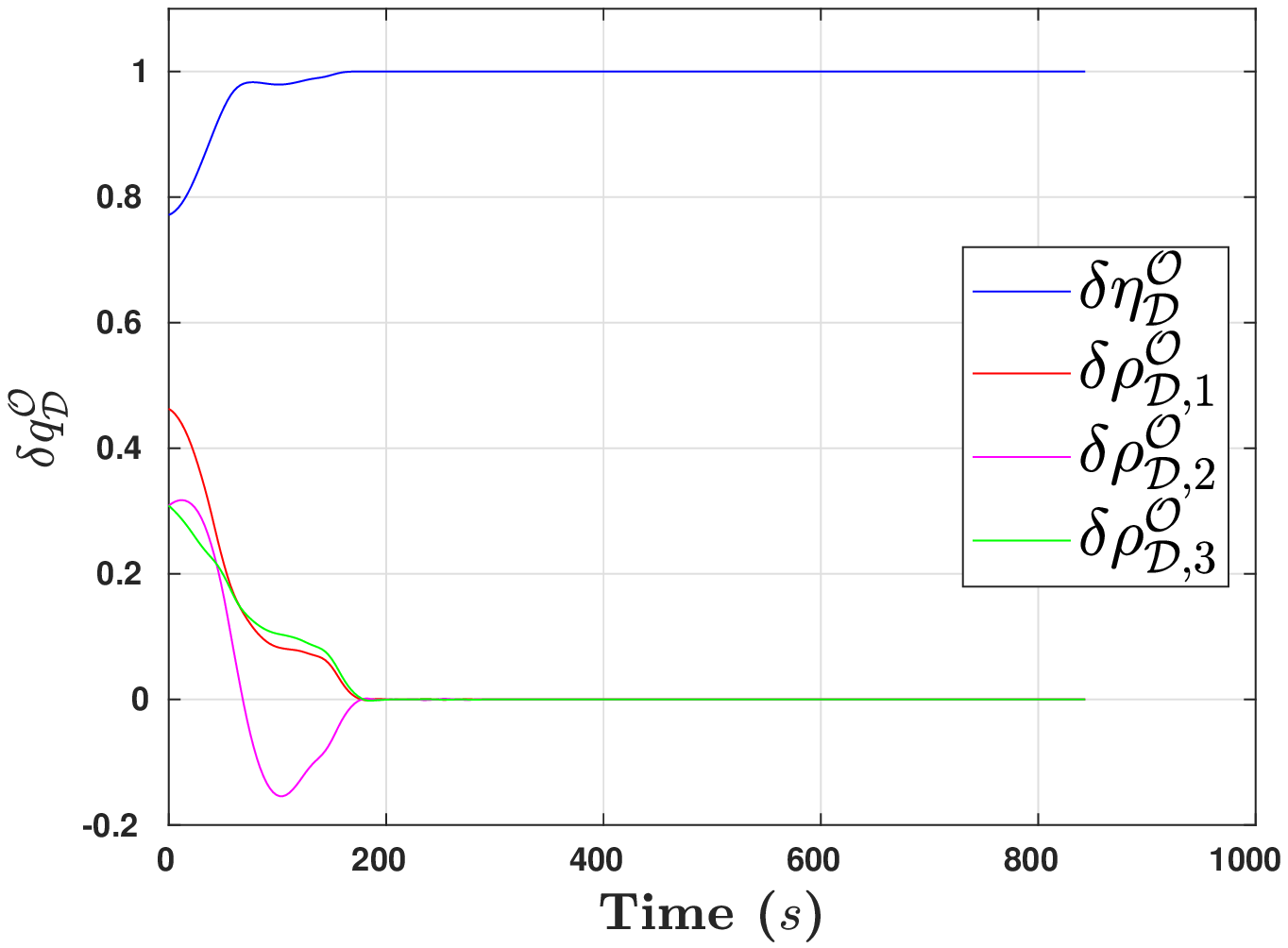}
            \caption[]%
            {{\small Error Quaternion,~$\delta q^{\coordframe{O}}_{\coordframe{D}}$}}    
            \label{fig:quaternion2}
        \end{subfigure}
        \hfill
        \begin{subfigure}[b]{0.24\textwidth}   
            \centering 
            \includegraphics[width=\textwidth]{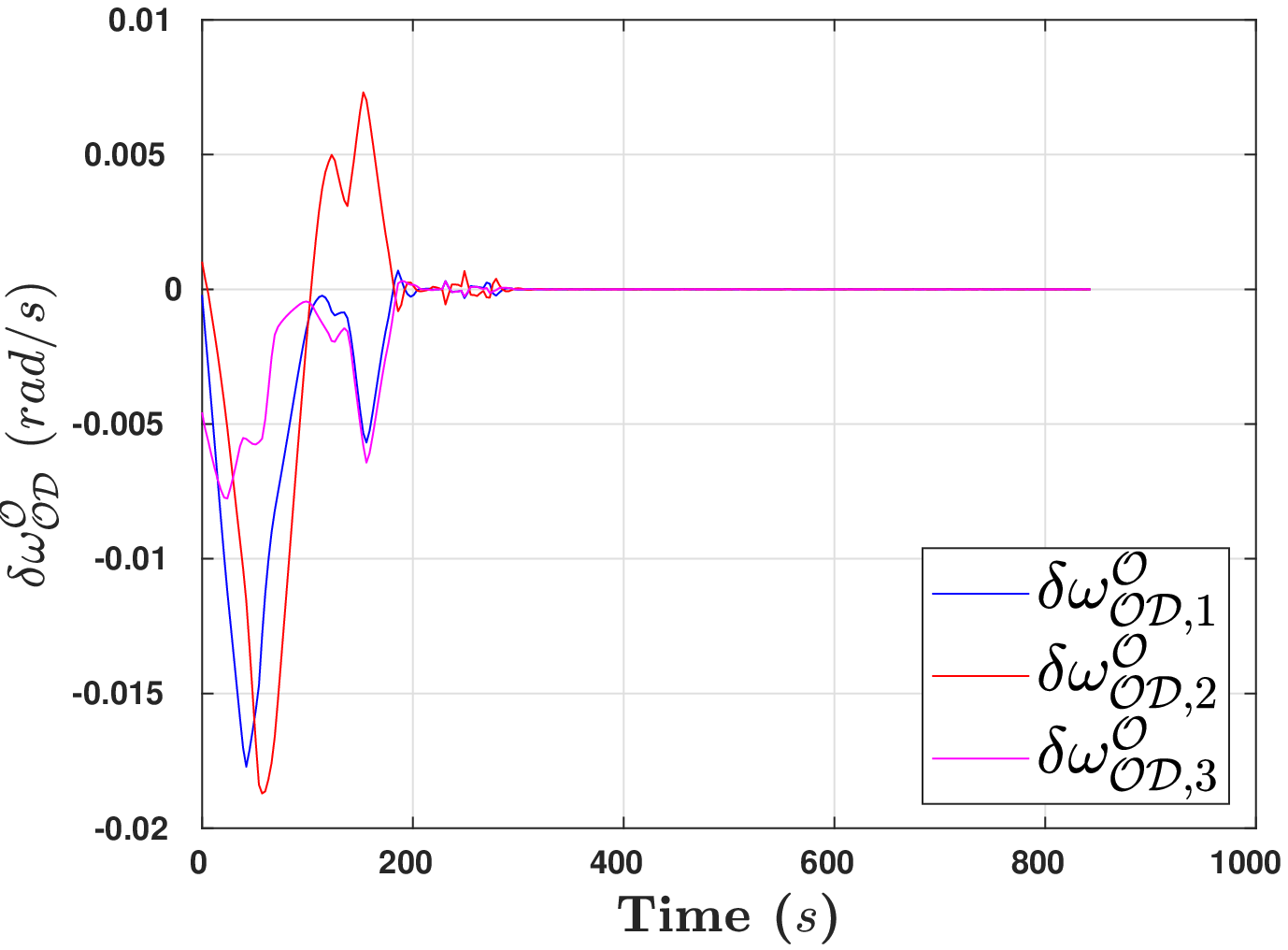}
            \caption[]%
            {{\small Error Angular Velocity,~$\delta\omega^{\coordframe{O}}_{\coordframe{O}\coordframe{D}}$}}    
            \label{fig:angularVel2}
        \end{subfigure}
        \caption[]
        {\small Plots of the perturbed attitudinal error, states, and torque inputs for $j_{max}=3$.} 
        \label{fig:attitude2}
\end{figure}

\begin{figure}
        \centering
        \begin{subfigure}[b]{0.23\textwidth}
            \centering
            \includegraphics[width=\textwidth]{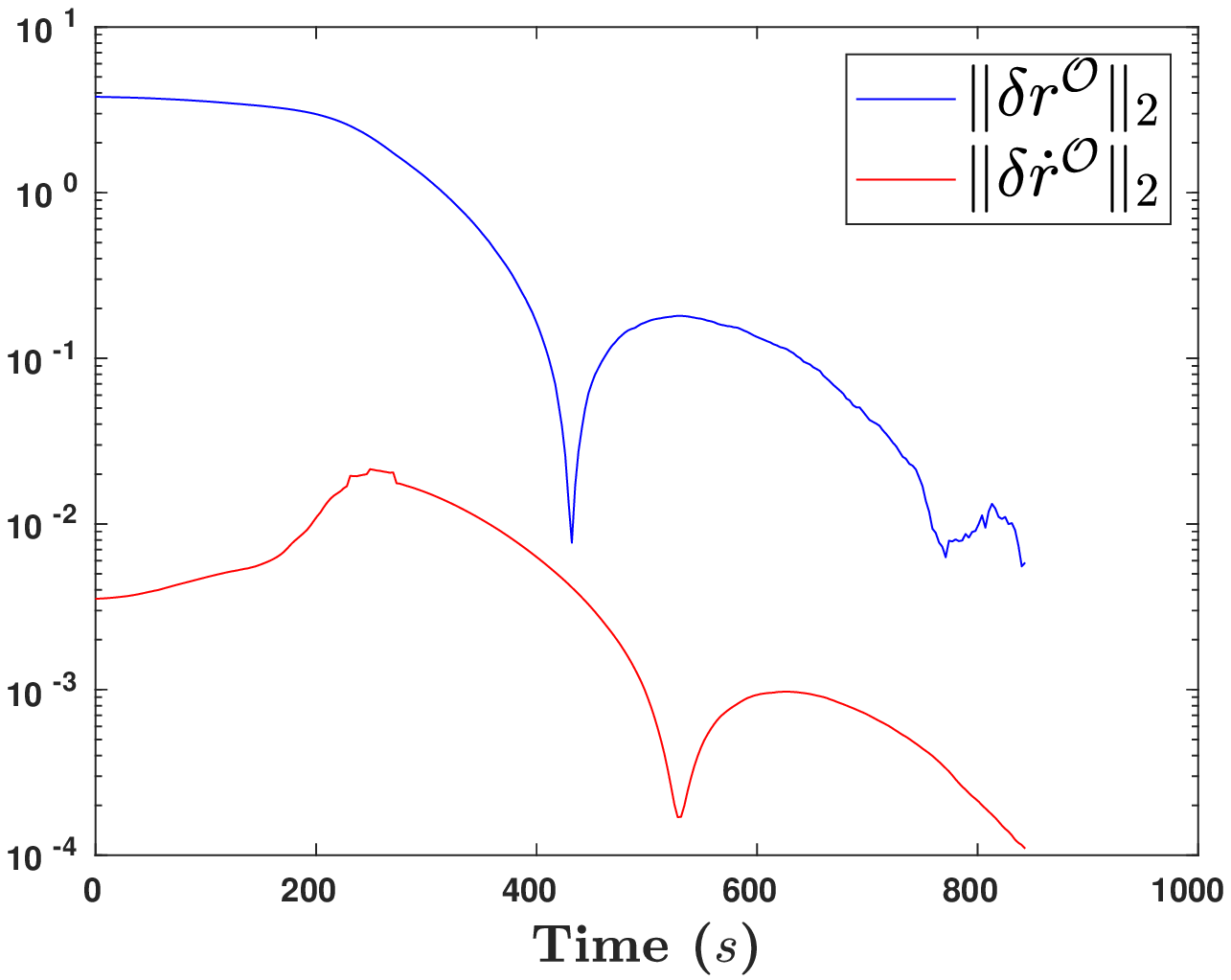}
            \caption[]%
            {{\small Translational Error}}    
            \label{fig:translationError2}
        \end{subfigure}
        \hfill
        \begin{subfigure}[b]{0.23\textwidth}  
            \centering 
            \includegraphics[width=\textwidth]{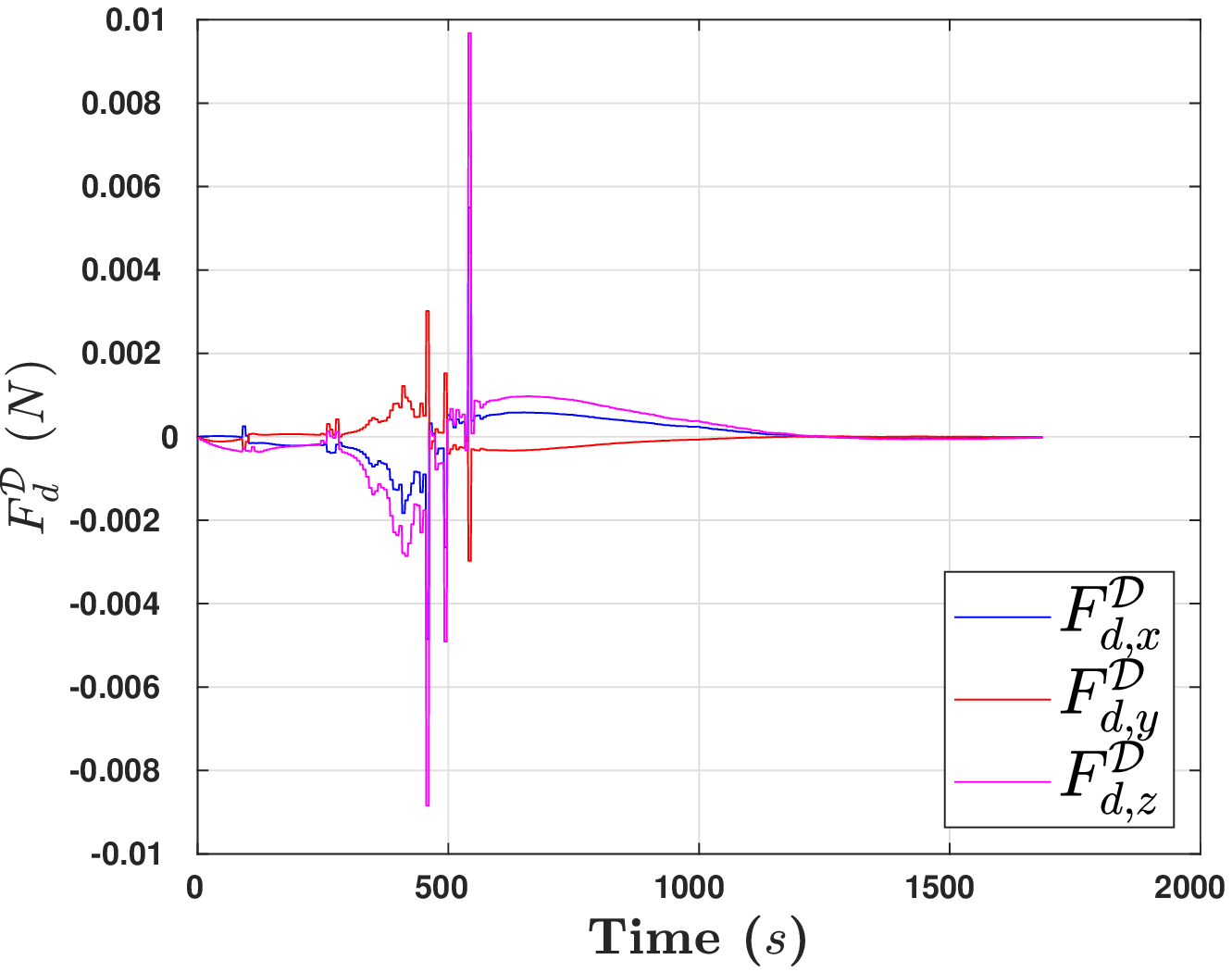}
            \caption[]%
            {{\small Thrust Inputs,~$F^{\coordframe{D}}_d \ (\text{N})$}}    
            \label{fig:thrust2}
        \end{subfigure}
        \vskip\baselineskip
        \begin{subfigure}[b]{0.23\textwidth}   
            \centering 
            \includegraphics[width=\textwidth]{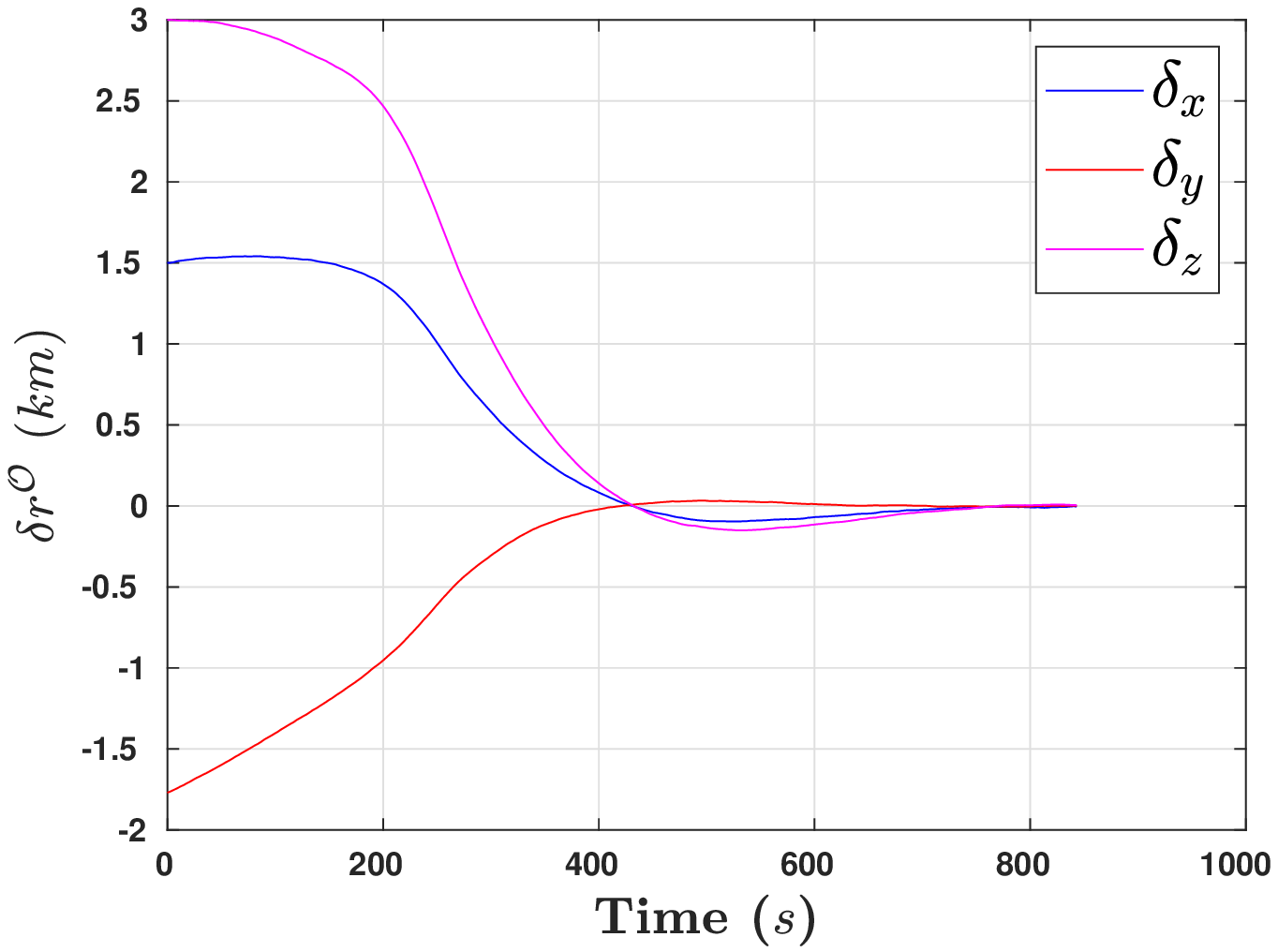}
            \caption[]%
            {{\small Position,~$\delta r^{\coordframe{O}} \ (\text{km})$}}    
            \label{fig:position2}
        \end{subfigure}
        \hfill
        \begin{subfigure}[b]{0.23\textwidth}   
            \centering 
            \includegraphics[width=\textwidth]{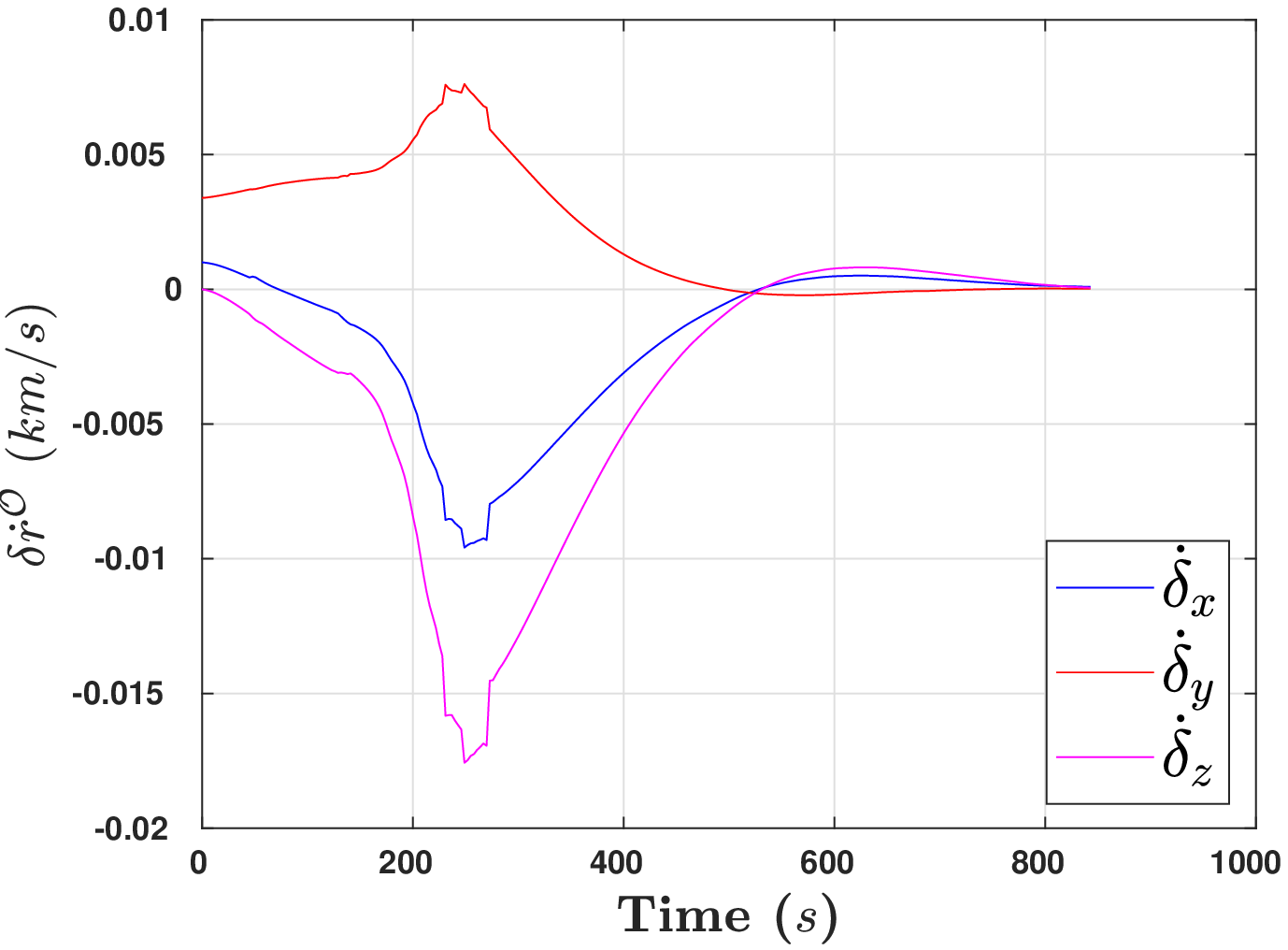}
            \caption[]%
            {{\small Linear Velocity,~$\delta \dot{r}^{\coordframe{O}} \ (\text{km/s})$}}    
            \label{fig:velocity2}
        \end{subfigure}
        \caption[]
        {\small Plots of the perturbed translational error, states, and thrust inputs when~$j_{max}=3$.} 
        \label{fig:translation2}
\end{figure}

The optimal error trajectory
takes longer to reach zero than the error trajectories
generated under the computational time constraint.
This differs from the unperturbed case and is perhaps counterintuitive,
but is easily explained. The optimal error trajectory takes
longer to reach zero because 
 the control inputs are computed using our dynamic model without knowledge of the perturbations. 
 Because future perturbations are unknown, exactly solving for an optimum
 while predicting future states does not necessarily give the best performance. 
 In Figures~\ref{fig:attitude2} and~\ref{fig:translation2} we show the states and inputs under perturbations for the initial conditions in Table~\ref{tb:parameters} and~$j_\max = 3$. 
 %\mh{Same as above: how many iterations would we need if we didn't constrain~$j_k$?}
 %\gb{same as above}
 
 These plots are similar to the unperturbed case but we can see the effect of perturbations, especially in the attitude error in Figure~\ref{fig:attitudeError2}.
 In all cases, the docking configuration is successfully reached. Table~\ref{tb:perturbTimes} shows that we achieve at least a 23.98\% reduction in average solving time per MPC loop for all values of~$j_\max$. Table~\ref{tb:perturbTimes} also shows a roughly 95\% reduction
 in the worst case solving time for all~$j_\max$.

\begin{comment}
    See~\citep[Section 8.6.2]{lavaei2007} for more detail on numerically bounding these perturbations.
\end{comment}

\begin{table} 
\caption{Solver times and reductions for Problem~\ref{prob:arpodMPC} 
with perturbations
}
\resizebox{\columnwidth}{10mm}{
\begin{tabular}{|| M{1cm} | M{2.2cm}| M{2.3cm} | M{2.2cm} | M{2.3cm} ||} 
 \hline
 $j_{\max}$  &   Average Time/Loop Reduction (s)&  Average Time/Loop Reduction (\%) & Maximum Loop Time Reduction (s) & Maximum Loop Time Reduction (\%)\\  [0.5ex] 
 \hline\hline
 2 & 0.3616 & 36.16\% & 26.64 & 95.80\% \\
 \hline
 3 & 0.2385  & 23.85\% & 26.33 & 94.70\%\\
 \hline
 4 & 0.2029 & 20.29\% & 26.66 & 95.86\% \\
 \hline
 5 & 0.2396 & 23.96\% & 26.65 & 95.86\% \\ [1ex] 
  \hline
\end{tabular}
}
\label{tb:perturbTimes}
\end{table}

\section{Conclusion} \label{sec:conclusion}
We proposed a time-constrained MPC control strategy for the 6 degree of freedom ARPOD problem, and showed that we can achieve the desired docking configuration under
restrictive computational time constraints. 
%Specifically, we showed that, even with the number of allowable computations
%reduced by more than~$90$\%, the docking configuration is still reached. 
Moreover, we demonstrated the robustness of 
docking using
time-constrained MPC where the states of the system were subject to perturbations.

\bibliography{ifacconf}             % bib file to produce the bibliography

\begin{thebibliography}{37}
\providecommand{\natexlab}[1]{#1}
\providecommand{\url}[1]{\texttt{#1}}
\providecommand{\urlprefix}{URL }
\expandafter\ifx\csname urlstyle\endcsname\relax
  \providecommand{\doi}[1]{doi:\discretionary{}{}{}#1}\else
  \providecommand{\doi}{doi:\discretionary{}{}{}\begingroup
  \urlstyle{rm}\Url}\fi

\bibitem[{ESA(2022)}]{ESAreport}
 (2022).
\newblock Esa’s annual space environment report.
\newblock Technical Report GEN-DB-LOG-00288-OPS-SD, ESA Space Debris Office,
  Darmstadt, Germany.

\bibitem[{Bernhard et~al.(2020)Bernhard, Choi, Rahmani, Chung, and
  Hadaegh}]{bernhard2020}
Bernhard, B., Choi, C., Rahmani, A., Chung, S.J., and Hadaegh, F. (2020).
\newblock Coordinated motion planning for on-orbit satellite inspection using a
  swarm of small-spacecraft.
\newblock In \emph{2020 IEEE Aerospace Conference}, 1--13.

\bibitem[{Bourdarie and Xapsos(2008)}]{bourdarie2008near}
Bourdarie, S. and Xapsos, M. (2008).
\newblock The near-earth space radiation environment.
\newblock \emph{IEEE transactions on nuclear science}, 55(4), 1810--1832.

\bibitem[{Cheng et~al.(2018)Cheng, Rivkin, Michel, Atchison, Barnouin, Benner,
  Chabot, Ernst, Fahnestock, Kueppers et~al.}]{cheng2018}
Cheng, A.F., Rivkin, A.S., Michel, P., Atchison, J., Barnouin, O., Benner, L.,
  Chabot, N.L., Ernst, C., Fahnestock, E.G., Kueppers, M., et~al. (2018).
\newblock Aida dart asteroid deflection test: Planetary defense and science
  objectives.
\newblock \emph{Planetary and Space Science}, 157, 104--115.

\bibitem[{Curtis(2013)}]{curtis2013orbital}
Curtis, H. (2013).
\newblock \emph{Orbital mechanics for engineering students}.
\newblock Butterworth-Heinemann.

\bibitem[{Di~Cairano et~al.(2012)Di~Cairano, Park, and Kolmanovsky}]{di2012}
Di~Cairano, S., Park, H., and Kolmanovsky, I. (2012).
\newblock Model predictive control approach for guidance of spacecraft
  rendezvous and proximity maneuvering.
\newblock \emph{Int. J. of Robust and Nonlinear Control}, 22(12), 1398--1427.

\bibitem[{Dong et~al.(2018)Dong, Hu, and Akella}]{dong2018}
Dong, H., Hu, Q., and Akella, M. (2018).
\newblock Dual-quaternion-based spacecraft autonomous rendezvous and docking
  under six-degree-of-freedom motion constraints.
\newblock \emph{J. of Guidance, Control, and Dynamics}, 41(5), 1150--1162.

\bibitem[{Falcone et~al.(2007)Falcone, Tufo, Borrelli, Asgari, and
  Tseng}]{falcone2007}
Falcone, P., Tufo, M., Borrelli, F., Asgari, J., and Tseng, H.E. (2007).
\newblock A linear time varying model predictive control approach to the
  integrated vehicle dynamics control problem in autonomous systems.
\newblock In \emph{2007 46th IEEE Conference on Decision and Control},
  2980--2985.

\bibitem[{Graichen and K{\"a}pernick(2012)}]{graichen2012}
Graichen, K. and K{\"a}pernick, B. (2012).
\newblock \emph{A real-time gradient method for nonlinear model predictive
  control}.
\newblock INTECH Open Access Publisher London.

\bibitem[{Graichen and Kugi(2010)}]{graichen2010}
Graichen, K. and Kugi, A. (2010).
\newblock Stability and incremental improvement of suboptimal mpc without
  terminal constraints.
\newblock \emph{IEEE Trans. on Automatic Control}, 55(11), 2576--2580.

\bibitem[{Gr{\"u}ne and Pannek(2017)}]{grune2017}
Gr{\"u}ne, L. and Pannek, J. (2017).
\newblock Nonlinear model predictive control.
\newblock In \emph{Nonlinear model predictive control}, 45--69. Springer.

\bibitem[{Hartley(2015)}]{hartley2015}
Hartley, E.N. (2015).
\newblock A tutorial on model predictive control for spacecraft rendezvous.
\newblock In \emph{2015 European Control Conference (ECC)}, 1355--1361. IEEE.

\bibitem[{Hartley and Maciejowski(2013)}]{hartley2013}
Hartley, E.N. and Maciejowski, J.M. (2013).
\newblock Graphical fpga design for a predictive controller with application to
  spacecraft rendezvous.
\newblock In \emph{52nd IEEE Conference on Decision and Control}, 1971--1976.

\bibitem[{Hogan and Schaub(2014)}]{hogan2014attitude}
Hogan, E.A. and Schaub, H. (2014).
\newblock Attitude parameter inspired relative motion descriptions for relative
  orbital motion control.
\newblock \emph{Journal of Guidance, Control, and Dynamics}, 37(3), 741--749.

\bibitem[{Hu et~al.(2021)Hu, Yang, Dong, and Zhao}]{hu2021}
Hu, Q., Yang, H., Dong, H., and Zhao, X. (2021).
\newblock Learning-based 6-dof control for autonomous proximity operations
  under motion constraints.
\newblock \emph{IEEE Trans. on Aerospace and Electronic Systems}, 57(6),
  4097--4109.

\bibitem[{Jewison(2017)}]{jewison2017guidance}
Jewison, C.M. (2017).
\newblock \emph{Guidance and Control for Multi-Stage Rendezvous and Docking
  Operations in the Presence of Uncertainty}.
\newblock Ph.D. thesis, Massachusetts Institute of Technology.

\bibitem[{Kuipers(1999)}]{kuipers1999quaternions}
Kuipers, J.B. (1999).
\newblock \emph{Quaternions and rotation sequences: a primer with applications
  to orbits, aerospace, and virtual reality}.
\newblock Princeton university press.

\bibitem[{Lee and Mesbahi(2017)}]{lee2017}
Lee, U. and Mesbahi, M. (2017).
\newblock Constrained autonomous precision landing via dual quaternions and
  model predictive control.
\newblock \emph{Journal of Guidance, Control, and Dynamics}, 40(2), 292--308.

\bibitem[{Leomanni et~al.(2014)Leomanni, Rogers, and Gabriel}]{leomanni2014}
Leomanni, M., Rogers, E., and Gabriel, S.B. (2014).
\newblock Explicit model predictive control approach for low-thrust spacecraft
  proximity operations.
\newblock \emph{Journal of Guidance, Control, and Dynamics}, 37(6), 1780--1790.

\bibitem[{Li et~al.(2017)Li, Yuan, Zhang, and Gao}]{li2017}
Li, Q., Yuan, J., Zhang, B., and Gao, C. (2017).
\newblock Model predictive control for autonomous rendezvous and docking with a
  tumbling target.
\newblock \emph{Aerospace Science and Technology}, 69, 700--711.

\bibitem[{Lovelly and George(2017)}]{lovelly2017comparative}
Lovelly, T.M. and George, A.D. (2017).
\newblock Comparative analysis of present and future space-grade processors
  with device metrics.
\newblock \emph{Journal of Aerospace Information Systems}, 14(3), 184--197.

\bibitem[{Lovelly(2017)}]{lovelly2017comparative_dissertation}
Lovelly, T.M. (2017).
\newblock \emph{Comparative Analysis of Space-Grade Processors}.
\newblock Ph.D. thesis, University of Florida.

\bibitem[{Malyuta et~al.(2021)Malyuta, Yu, Elango, and
  A{\c{c}}{\i}kme{\c{s}}e}]{malyuta2021}
Malyuta, D., Yu, Y., Elango, P., and A{\c{c}}{\i}kme{\c{s}}e, B. (2021).
\newblock Advances in trajectory optimization for space vehicle control.
\newblock \emph{Annual Reviews in Control}, 52, 282--315.

\bibitem[{Mayne(2014)}]{mayne2014}
Mayne, D.Q. (2014).
\newblock Model predictive control: Recent developments and future promise.
\newblock \emph{Automatica}, 50(12), 2967--2986.

\bibitem[{Ogilvie et~al.(2008)Ogilvie, Allport, Hannah, and
  Lymer}]{ogilvie2008}
Ogilvie, A., Allport, J., Hannah, M., and Lymer, J. (2008).
\newblock Autonomous robotic operations for on-orbit satellite servicing.
\newblock In \emph{Sensors and Systems for Space Applications II}, volume 6958,
  50--61.

\bibitem[{Pavlov et~al.(2019)Pavlov, Shames, and Manzie}]{pavlov2019}
Pavlov, A., Shames, I., and Manzie, C. (2019).
\newblock Early termination of nmpc interior point solvers: Relating the
  duality gap to stability.
\newblock In \emph{2019 18th European Control Conference (ECC)}, 805--810.
  IEEE.

\bibitem[{Petersen et~al.(2014)Petersen, Jaunzemis, Baldwin, Holzinger, and
  Kolmanovsky}]{petersen2014}
Petersen, C., Jaunzemis, A., Baldwin, M., Holzinger, M., and Kolmanovsky, I.
  (2014).
\newblock Model predictive control and extended command governor for improving
  robustness of relative motion guidance and control.
\newblock In \emph{Proc. AAS/AIAA space flight mechanics meeting}.

\bibitem[{Richards and How(2003)}]{richards2003}
Richards, A. and How, J.P. (2003).
\newblock Model predictive control of vehicle maneuvers with guaranteed
  completion time and robust feasibility.
\newblock In \emph{2003 American Control Conference}, volume~5, 4034--4040.

\bibitem[{Soderlund and Phillips(2022)}]{soderlund2022switching}
Soderlund, A.A. and Phillips, S. (2022).
\newblock Autonomous rendezvous and proximity operations of an underactuated
  spacecraft via switching controls.
\newblock In \emph{AIAA SCITECH 2022 Forum}, 0956.

\bibitem[{Soderlund et~al.(2021)Soderlund, Phillips, Zaman, and
  Petersen}]{soderlund2021autonomous}
Soderlund, A.A., Phillips, S., Zaman, A., and Petersen, C.D. (2021).
\newblock Autonomous satellite rendezvous and proximity operations via
  geometric control methods.
\newblock In \emph{AIAA Scitech 2021 Forum}, 0075.

\bibitem[{Sun and Huo(2015)}]{sun2015}
Sun, L. and Huo, W. (2015).
\newblock 6-dof integrated adaptive backstepping control for spacecraft
  proximity operations.
\newblock \emph{IEEE Transactions on Aerospace and Electronic Systems}, 51(3),
  2433--2443.

\bibitem[{Trivailo et~al.(2009)Trivailo, Wang, and Zhang}]{trivailo2009}
Trivailo, P.M., Wang, F., and Zhang, H. (2009).
\newblock Optimal attitude control of an accompanying satellite rotating around
  the space station.
\newblock \emph{Acta Astronautica}, 64(2-3), 89--94.

\bibitem[{W{\"a}chter and Biegler(2006)}]{wachter2006}
W{\"a}chter, A. and Biegler, L.T. (2006).
\newblock On the implementation of an interior-point filter line-search
  algorithm for large-scale nonlinear programming.
\newblock \emph{Mathematical programming}, 106(1), 25--57.

\bibitem[{Wang et~al.(2018)Wang, Wang, and Zhang}]{wang2018}
Wang, X., Wang, Z., and Zhang, Y. (2018).
\newblock Model predictive control to autonomously approach a failed
  spacecraft.
\newblock \emph{Int. J. of Aerospace Engineering}, 2018.

\bibitem[{Wang and Ji(2019)}]{wang2019}
Wang, Y. and Ji, H. (2019).
\newblock Integrated relative position and attitude control for spacecraft
  rendezvous with iss and finite-time convergence.
\newblock \emph{Aerospace Science and Technology}, 85, 234--245.

\bibitem[{Yang and Stoll(2019)}]{yang2019}
Yang, J. and Stoll, E. (2019).
\newblock Adaptive sliding mode control for spacecraft proximity operations
  based on dual quaternions.
\newblock \emph{Journal of Guidance, Control, and Dynamics}, 42(11),
  2356--2368.

\bibitem[{Zhou et~al.(2020)Zhou, Liu, and Cai}]{zhou2020}
Zhou, B.Z., Liu, X.F., and Cai, G.P. (2020).
\newblock Motion-planning and pose-tracking based rendezvous and docking with a
  tumbling target.
\newblock \emph{Advances in Space Research}, 65(4), 1139--1157.

\end{thebibliography}
                                                     % with bibtex (preferred)

                                                                         % in the appendices.
\end{document}